\documentclass[12pt]{amsart}

\usepackage{cite}
\usepackage{graphicx}
\usepackage{amssymb}

\newcommand{\0}{\mathbf{0}}
\newcommand{\1}{\mathbf{1}}
\newcommand{\cle}{\preceq}
\newcommand{\opl}{{\oplus}}
\newcommand{\bset}{\mathbf{B}}
\newcommand{\cset}{\mathbf{C}}
\newcommand{\pset}{\mathbf{P}}
\newcommand{\rset}{\mathbf{R}}
\newcommand{\maA}{\mathcal{A}}
\newcommand{\CF}{\mathcal{F}}
\newcommand{\rmin}{\mathbf{R}_{\min}}
\newcommand{\rmax}{\mathbf{R}_{\max}}
\newcommand{\Log}{\mathop{\mathrm{Log}}}
\newcommand{\suplim}{\sup\limits}
\newcommand{\sumlim}{\sum\limits}
\newcommand{\maxlim}{\max\limits}
\newcommand{\pd}[2]{\dfrac{\partial#1}{\partial#2}}
\newcommand{\CalB}{\mathcal{B}}

\renewcommand{\Re}{\mathop{\text{Re}}}

\newtheorem*{thm}{Theorem}
\newtheorem{prop}{Proposition}
\newtheorem{cor}{Corollary}

\newcommand{\maF}{\mathcal{F}}
\newcommand{\maD}{\mathcal{D}}

\begin{document}

\title[The Maslov dequantization]
{The Maslov Dequantization, Idempotent and Tropical Mathematics:
a Brief Introduction}

\author{G.~L.~Litvinov}
\address{Independent University of Moscow\\ 
Bol'shoi Vlasievskii per., 11\\ 
Moscow 121002, Russia }
\email{islc@dol.ru} 

\dedicatory{To Anatoly Vershik with admiration and gratitude}

\date{}

\thanks{This work has been supported by the RFBR grant 05--01--00824
and the Erwin Schr{\"o}\-din\-ger International Institute for
Mathematical Physics (ESI).\\[2ex]
To be published in \textit{Journal of Mathematical Sciences}}

\begin{abstract}
This paper is a brief introduction to idempotent and tropical mathematics.
Tropical mathematics can be treated as a result of the so-called Maslov
dequantization of the traditional mathematics over numerical fields as the
Planck constant $\hbar$ tends to zero taking imaginary values.
\end{abstract}

\maketitle

This paper is a brief introduction, practically without exact theorems and
proofs, to the Maslov dequantization and idempotent and tropical
mathematics.  Our list of references is not complete (not at all).
Additional references can be found, e.g., in the
electronic archive 
\[
  \text{\texttt{http://arXiv.org}}
\]
and in \cite{BaCoOlQu92,Bu94,%
Ca79,CoGaQu99,CoQu94,Cu79,Cu95,CuMe80,DeDo98,DeDo2001,%
Fl2004,Gla2002,Gol99,Gol2000,Gol2003,GoMi79,GoMi2001,Gun98a,%
Gun98b,KiRo2004,Ki2001,KlMePa2000,Ko2001,KoMa97,LiMa95,LiMa98,%
LiMaSh2001,LiSo2001,Pap2002,ZiU81}.
The present brief survey is an extended version of the paper
\cite{Li2005}.

The author thanks M.\ Akian, Ya.\ I.\ Belopolskaya, P.\ Butkovi{\v c}, G.\ Cohen, S.\ Gaubert,
R.\ A.\ Cuninghame-Green, P.\ Del Moral, W.\ H.\ Fleming,
J.\ S.\ Golan, M.\ Gondran, I.\ Itenberg, Y.\ Katsov, V.\ N.\ Kolokoltsov,
P.\ Loreti, Yu. I.\ Manin, G.\ Mascari,
W.\ M.\ Mac\-Eneaney, G.\ Panina, E.\ Pap, M.\ Pedicini, H.\ Prade, A.\ A.\ Puhalskii, 
J.-P.\ Quadrat, M.\ A.\ Roytberg, G.\ B.\ Shpiz, I.\ Singer,
and O.\ Viro for useful contacts and references.  Special thanks are due
to V.\ P.\ Maslov for his crucial help and support
and to A.\ N.\ Sobolevski\u\i \ for his great help including, but not
limited to, three pictures presented below.

%{\bf 1. Some basic ideas.}
\section{Some basic ideas}

Idempotent mathematics is based on replacing the usual arithmetic
operations with a new set of basic operations (such as maximum
or minimum), that is on replacing numerical fields by idempotent
semirings and semifields. Typical examples are given by the
so-called max-plus algebra $\rmax$ and the min-plus algebra
$\rmin$. Let $\rset$ be the field of real numbers. Then
${\rmax}={\rset} \cup \{-\infty\}$ with operations $x\oplus y=\max
\{x,y\}$ and $x\odot y=x+y$. Similarly ${\rmin}={\rset}\cup
\{+\infty\}$ with the operations $\oplus=\min$, $\odot=+$. The new
addition $\oplus$ is idempotent, i.e., $x\oplus x=x$ for all
elements $x$.

Many authors (S.\ C.\ Kleene, S.\ N.\ N.\ Pandit, N.\ N.\ Vorobjev, 
B.\ A.\ Carre,
R.\ A.\ Cu\-ning\-ha\-me-Gre\-en, K.\ Zimmermann, U.\ Zimmermann, M.\ Gondran, 
F.\ L.\ Baccelli, G.\ Cohen, S.\ Gaubert,
G.\ J.\ Olsder, J.-P.\ Quadrat, and others) used idempotent semirings and
matrices over these semirings for solving some applied problems in
computer science and discrete mathematics, starting from the classical
paper of S.\ C.\ Kleene \cite{Kle56}.  The modern 
\emph{idempotent analysis} (or \emph{idempotent calculus}, or \emph{idempotent
mathematics}) was founded by V.\ P.\ Maslov and his collaborators in the 
1980s in Moscow; see,
e.g., \cite{Mas86,Mas87a,Mas87b,MaKo94,MaSa92,MaVo88,KoMa97}.
Some preliminary results are due to E. Hopf and G. Choquet, see
\cite{Cho55,Ho50}.

Idempotent mathematics can be treated as a result of a
dequantization of the traditional mathematics over numerical
fields as the Planck constant $\hbar$ tends to zero taking
imaginary values.  This point of view was presented by 
G.~L.~Litvinov and V.~P.~Maslov \cite{LiMa95,LiMa96,LiMa98}, see also
\cite{LiMaSh2001,LiMaSh2002}.  In other words,
idempotent mathematics is an asymptotic version of the 
traditional
 mathematics over the fields of real and complex numbers.

The basic paradigm is expressed in terms of an {\it idempotent
correspondence principle}.  This principle is closely related to the
well-known correspondence
principle of N.~Bohr in quantum theory.
Actually, there exists a heuristic correspondence between
important, interesting, and useful constructions and results of the
traditional mathematics over fields and analogous constructions and results
over idempotent semirings and semifields (i.e., semirings and semifields
with idempotent addition).

\begin{figure}[t]
  \label{fig:quant}
  \centering
  \includegraphics[width=10cm]{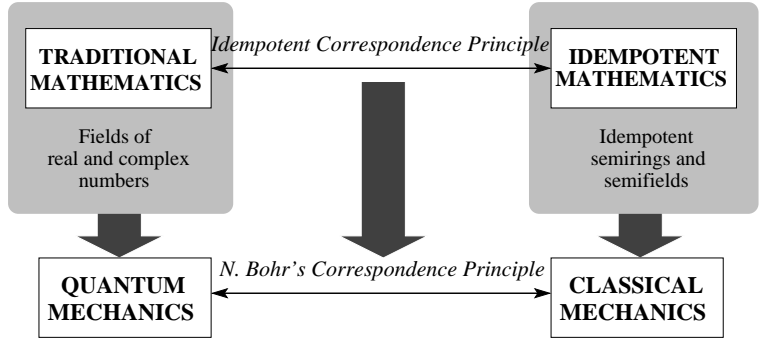}
  \caption{Relations between idempotent and traditional mathematics.}
\end{figure}

A systematic and consistent application of the idempotent
correspondence principle leads to a variety of results, often
quite unexpected.  As a result, in parallel with the traditional
mathematics over fields, its ``shadow,'' the idempotent
mathematics, appears.  This ``shadow'' stands approximately in the
same relation to the traditional mathematics as does classical physics
to quantum theory, see Fig.~\ref{fig:quant}. 
In many respects idempotent mathematics is
simpler than the traditional one.  However the transition from
traditional concepts and results to their idempotent analogs is
often nontrivial.

%%% {\bf 2. Semirings, semifields, and dequantization.} 
\section{Semirings, semifields, and dequantization}

Consider a
set $S$ equipped with two algebraic operations: {\it addition} $\oplus$
and {\it multiplication} $\odot$. It is a {\it semiring} if the following
conditions are satisfied:
\begin{itemize}
\item the addition $\oplus$ and the multiplication $\odot$ are
associative;
\item the addition $\oplus$ is commutative;
\item the multiplication $\odot$ is distributive with respect to
the addition $\oplus$: \[x\odot(y\oplus z)=(x\odot y)\oplus(x\odot z)
\quad\text{and}\quad
(x\oplus y)\odot z=(x\odot z)\oplus(y\odot z)\] 
for all $x,y,z\in S$.
\end{itemize}
A {\it unity} of a semiring $S$ is an element $\1\in S$ such that
$\1\odot x=x\odot\1=x$ for all $x\in S$. A {\it zero} of a semiring
$S$ is an element $\0\in S$ such that $\0\neq\1$ and $\0\oplus x=x$,
$\0\odot x=x\odot \0=\0$ for all $x\in S$. A semiring $S$ is called an
{\it idempotent semiring} if $x\oplus x=x$ for all $x\in S$. A
semiring $S$ with neutral elements $\0$ and $\1$ is called a {\it
semifield} if every nonzero element of $S$ is invertible.
Note that dio{\"\i}ds in the sense of \cite{BaCoOlQu92,GoMi79,%
GoMi2001}, quantales in the sense of \cite{Ro90,Ro96},
and inclines in the sense of \cite{KiRo2004} are examples of
idempotent semirings.

\begin{figure}[t]
  \centering
  \includegraphics[height=10cm]{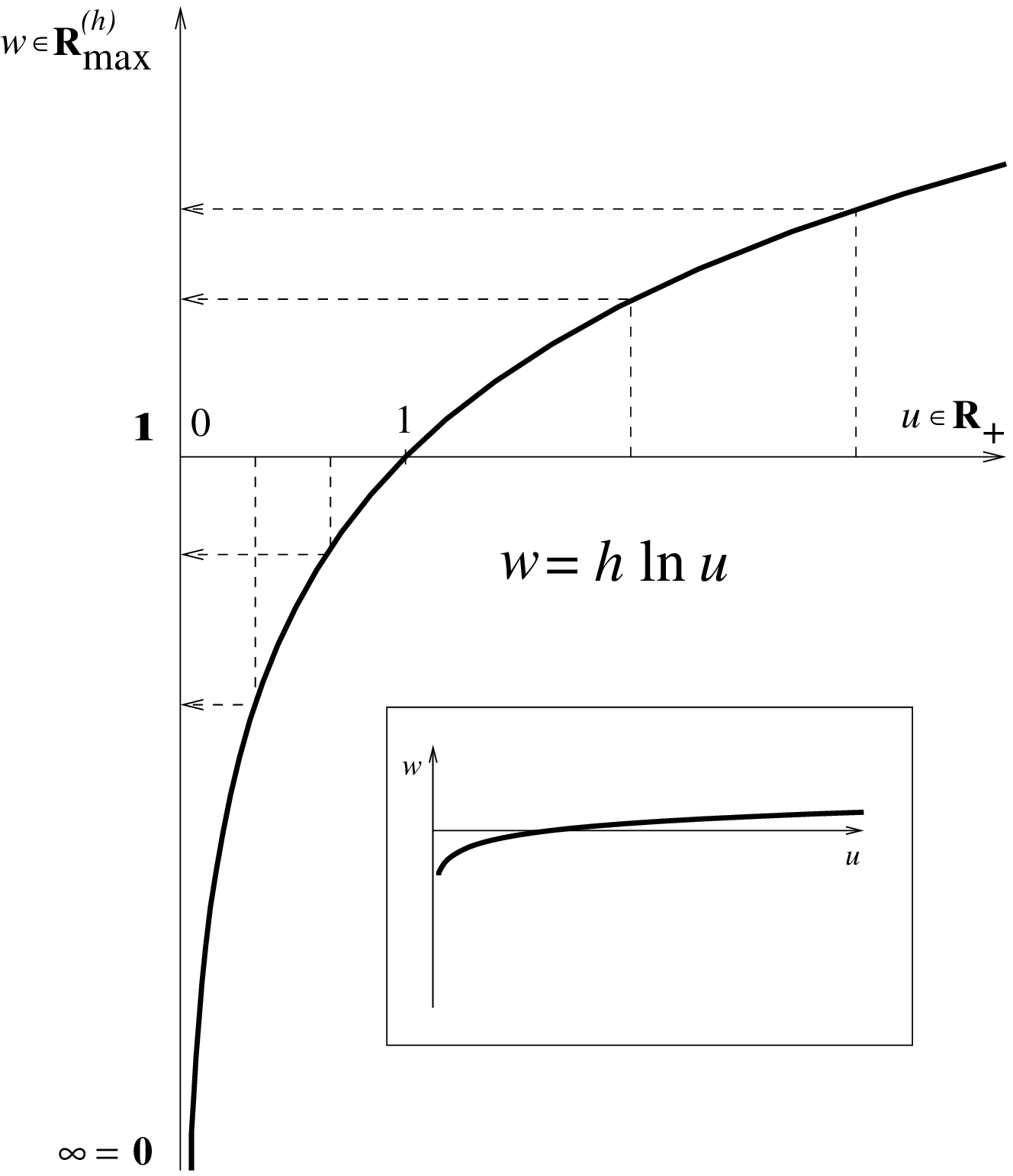}
  \label{fig:dequant}
  \caption{Deformation of~\(\rset_+\) to~\(\rset^{(h)}\).  Inset: same
  for a small value of~\(h\).}
\end{figure}

Let $\rset$ be the field of real numbers and $\rset_+$ the
semiring of all nonnegative real numbers (with respect to the
usual addition and multiplication). The change of variables $x
\mapsto u = h \ln x$, $h > 0$, defines a map $\Phi_h \colon
\rset_+ \to S = \rset \cup \{-\infty\}$.  Let the addition and
multiplication operations be mapped from $\rset$ to $S$ by
$\Phi_h$, i.e., let $u \oplus_h v = h \ln (\exp (u/h) +
\exp(v/h))$, $u \odot v = u + v$, $\0 = -\infty = \Phi_h(0)$,
$\1 = 0 = \Phi_h(1)$, so $S$ has a structure of a semiring
$\rset^{(h)}$ isomorphic to $\rset_+$; see Fig.~\ref{fig:dequant}.

It can easily be checked that $u \oplus_h
v \to \max \{u,v\}$ as $h \to 0$ and that $S$ forms a semiring with
respect to addition $u \oplus v = \max \{u,v\}$ and multiplication
$u \odot v = u + v$ with zero $\0 = -\infty$ and unit $\1 = 0$.
Denote this semiring by $\rmax$; it is {\it idempotent}, i.e., $u
\oplus u = u$ for all its elements. The semiring $\rset_{\max}$ is
actually a semifield. The analogy with quantization is obvious;
the parameter $h$ plays the r\^{o}le of the Planck constant, so
$\rset_+$ can be viewed as a ``quantum object'' and
$\rmax$ as the result of its ``dequantization.'' A similar
procedure  (for $h<0$) gives the semiring 
$\rmin = \rset \cup \{+\infty\}$ with
the operations $\oplus = \min$, $\odot = +$; in this case $\0 =
+\infty$, $\1$ = 0. The semirings $\rmax$ and $\rset_{\min}$ are
isomorphic.  This passage to $\rmax$ or $\rmin$ is called the
{\it Maslov dequantization}. It is clear that the corresponding
passage from $\cset$ or $\rset$ to $\rset_{\max}$ is generated
by the Maslov dequantization and the map $x\mapsto |x|$.
By misuse of language, {\it we shall also call this passage the
Maslov dequantization}. Connections with physics and
the meaning of imaginary values of the Planck constant are discussed 
below (Section 6) and in
\cite{LiMaSh2001,LiMaSh2002}. The idempotent semiring $\rset \cup
\{-\infty\} \cup \{+\infty\}$ with the operations $\oplus = \max$,
$\odot = \min$ can be obtained as a result of a ``second
dequantization'' of $\cset$, $\rset$ or $\rset_+$. Dozens of interesting
examples of nonisomorphic idempotent semirings may be cited as
well as a number of standard methods of deriving new semirings
from these (see, e.g., \cite{CoGaQu2004,Gol99,Gol2003,%
GoMi79,GoMi2001,Gun98a,Gun98b,%
LiMa98,LiMaSh2001}). The so-called
{\it idempotent dequantization} is a generalization of the Maslov
dequantization; this is a passage from fields to idempotent
semifields and semirings in mathematical constructions and
results.

The Maslov dequantization is related to the well-known logarithmic
transformation that was used, e.g., in the classical papers of
E.~Schr{\"o}\-din\-ger \cite{Sch26} and E.~Hopf \cite{Ho50}. 
The term `Cole-Hopf
transformation' is also used. The subsequent progress of E.~Hopf's
ideas has culminated in the well-known vanishing viscosity method
and the method of viscosity solutions, see, e.g., \cite{BaCa97,%
CaLi97,FlSo93,MaKo94,McC2005,Rou2005,Su95}.

%%%{\bf 3. Tropical semirings and tropical mathematics.}
\section{Terminology: Tropical semirings and tropical mathematics}

The term `tropical semirings' was introduced in computer
science to denote discrete versions of the max-plus algebra
$\rmax$ or min-plus algebra $\rmin$ 
and their subalgebras; (discrete) semirings of this type
were called tropical semirings by Dominic Perrin in honour of Imre
Simon (who is a Brasilian mathematician and computer
scientist) because of his pioneering activity in this area,
see \cite{Pin98}.

More recently the situation and terminology have changed. For the
most part of modern authors `tropical' means `over $\rset_{\max}$ 
(or $\rset_{\min}$)' and tropical semirings are idempotent
semifields $\rset_{\max}$ and $\rset_{\min}$. The 
terms `max-plus' and `min-plus' are often used in the same
sense. Now the term `tropical mathematics' usually means
`mathematics over $\rset_{\max}$ or $\rset_{\min}$', see,
e.g., \cite{Ale2005,BeFoZe96,BeZe2001,BlYu2005,Dev2004,DeSaSt2004,%
FoGo2003a,FoGo2003b,%
GaMa2005a,GaMa2005b,ItKhSh2003,ItKhSh2003a,ItKhSh2004,ItKhSh2005,Izh2005a,%
Izh2005b,Jo2004,KiRo2005,KiRo2005a,Mi2001,Mi2003,Mi2004,Mi2005,NiSie2004,%
NoYa2002,PaSt2003,PaSt2004,RiStTh2005,Shu2002,Shu2004,So2005,Sp2004,SpSt2003,%
SpSt2004,SpWi2003,Te2004,The2004,Vi2004}.
Terms
`tropicalization' and `tropification' (see, e.g., \cite{Ki2001})
mean exactly dequantization and quantization in our sense. In
any case, tropical mathematics is a natural and very important
part of idempotent mathematics. Some well known constructions and
results of idempotent mathematics were repeated in the framework of
tropical mathematics (and especially tropical linear algebra).

Note that in papers \cite{Vor63,Vor67,Vor70} N.~N.~Vorobjev developed
a version of idempotent linear algebra (with important 
applications, e.g., to mathematical economics) and predicted 
many aspects of the future extended theory.
He used the terms `extremal algebras'
and `extremal mathematics' for idempotent semirings and 
idempotent mathematics.  Unfortunately, 
N.~N.~Vorobjev's papers and ideas were forgotten for a long period, so
his remarkable terminology is not in use any more.

%%% {\bf 4. Idempotent algebra and linear algebra.} 
\section{Idempotent algebra and linear algebra}

The first
known paper on idempotent linear algebra is due to S.~Klee\-ne
\cite{Kle56}. Systems of linear algebraic equations over an
exotic idempotent semiring of all formal languages over a
fixed finite alphabet are examined in this work; however, 
S.~Kleene's ideas are very general and universal.  Since then, dozens of
authors investigated matrices with coefficients belonging to
an idempotent semiring and the corresponding applications to
discrete mathematics, computer science, computer languages,
linguistic problems, finite automata, optimization problems
on graphs, discrete event systems and Petry nets, stochastic
systems, computer performance evaluation, computational problems
etc. This subject is very well known and well presented in the 
corresponding literature, see, e.g., \cite{BaCoOlQu92,Bu94,%
Bu2005,Ca79,CoGaQu99,CoQu94,Cu79,Cu95,DuSa92,Gla2002,%
Gol99,Gol2000,Gol2003,GoMi79,GoMi2001,Gun98a,Gun98b,%
KiRo2004,KoMa97,Ko2001,LiMa95,LiMa98,LiMa2005,LiMaE2000,%
LiMaRo2000,LiSo2000,LiSo2001,%
MaKo94,Vor63,Vor67,Vor70,Wa2005,ZiU81}.

Idempotent abstract algebra
is not so well developed yet (on the other hand, from a formal
point of view, the lattice theory and the theory of ordered
groups and semigroups are parts of idempotent algebra). However,
there are many interesting results and applications presented, e.g., 
in \cite{Cu79,Cu95,CuMe80,Ka2004,Ro90,Ro96,Shp2000}.

In particular, an idempotent version of the main theorem of
algebra holds \cite{CuMe80,Shp2000} for radicable idempotent
semifields (a semiring $A$ is {\it radicable} if the equation
$x^n = a$ has a solution $x\in A$ for any $a\in A$ and any
positive integer $n$). It is proved that $\rset_{\max}$ and other
radicable semifields are algebraically closed in a natural sense
\cite{Shp2000}.

Over the last years, tropical algebraic geometry attracted
a lot of attention. 
This subject will be briefly discussed below (Section 11).

%%%{\bf 5. Idempotent analysis.} 
\section{Idempotent analysis}

Idempotent analysis was initially
constructed by V.~P.~Maslov and his collaborators and then 
developed by many authors. The subject is presented in the book of
V.~N.~Kolokoltsov and V.~P.~Maslov \cite{KoMa97} (a version of
this book in Russian \cite{MaKo94} was published in 1994).

Let $S$ be an arbitrary semiring with idempotent addition $\oplus$
(which is always assumed to be commutative), multiplication
$\odot$, zero $\0$, and unit $\1$. The set $S$ is supplied with
the {\it standard partial order\/}~$\cle$: by definition, $a \cle
b$ if and only if $a \oplus b = b$. Thus all elements of $S$ are
nonnegative: $\0 \cle$ $a$ for all $a \in S$. Due to the existence
of this order, idempotent analysis is closely related to the lattice
theory, theory of vector lattices, and theory of ordered
spaces. Moreover, this partial order allows to model a number of
basic ``topological'' concepts and results of idempotent analysis at the purely
algebraic level; this line of reasoning was examined
systematically in \cite{LiMaSh98,LiMaSh99,LiMaSh2001,LiMaSh2002,%
LiSh2002} and \cite{CoGaQu2004}.

Calculus deals mainly with functions whose values are numbers. The
idempotent analog of a numerical function is a map $X \to S$,
where $X$ is an arbitrary set and $S$ is an idempotent semiring.
Functions with values in $S$ can be added, multiplied by each
other, and multiplied by elements of $S$ pointwise.

The idempotent analog of a linear functional space is a set of $S$-valued
functions that is closed under addition of functions and multiplication of
functions by elements of $S$, or an $S$-semimodule. Consider, e.g., the
$S$-semimodule $B(X, S)$ of all functions $X \to S$ that are bounded in
the sense of the standard order on $S$.

If $S = \rmax$, then the idempotent analog of integration is defined by the
formula
$$
I(\varphi) = \int_X^{\oplus} \varphi (x)\, dx     = \sup_{x\in X}
\varphi (x),\eqno{(1)}
$$
where $\varphi \in B(X, S)$. Indeed, a Riemann sum of the form
$\sumlim_i \varphi(x_i) \cdot \sigma_i$ corresponds to the expression
$\bigoplus\limits_i \varphi(x_i) \odot \sigma_i = \maxlim_i \{\varphi(x_i)
+ \sigma_i\}$, which tends to the right-hand side of~(1) as $\sigma_i \to
0$. Of course, this is a purely heuristic argument.

Formula~(1) defines the \emph{idempotent} (or \emph{Maslov}) \emph{integral} not only for 
functions taking
values in $\rmax$, but also in the general case when any of bounded
(from above) subsets of~$S$ has the least upper bound.

An \emph{idempotent} (or \emph{Maslov}) \emph{measure} on $X$ is defined by 
$m_{\psi}(Y) = \suplim_{x \in Y}
\psi(x)$, where $\psi \in B(X,S)$, $Y\subseteq X$. The integral with respect to this
measure is defined by
$$
   I_{\psi}(\varphi)
    = \int^{\oplus}_X \varphi(x)\, dm_{\psi}
    = \int_X^{\oplus} \varphi(x) \odot \psi(x)\, dx
    = \sup_{x\in X} (\varphi (x) \odot \psi(x)).
    \eqno{(2)}
$$

Obviously, if $S = \rmin$, then the standard order $\cle$ is opposite to
the conventional order $\le$, so in this case equation~(2)
assumes the form
$$
   \int^{\oplus}_X \varphi(x)\, dm_{\psi}
    = \int_X^{\oplus} \varphi(x) \odot \psi(x)\, dx
    = \inf_{x\in X} (\varphi (x) \odot \psi(x)),
$$
where $\inf$ is understood in the sense of the conventional order $\le$.

Note that the so-called pseudo-analysis (see, e.g., the survey paper of
E.~Pap \cite{Pap2005}) is related to a
special part of idempotent analysis; however, this
pseudo-analysis is not a proper part of idempotent 
mathematics in the general case.  Some generalizations of 
Maslov measures are discussed in \cite{KlMePa2000,Pap2002}.

%% {\bf 6. The superposition principle and linear problems.} 
\section{The superposition principle and linear problems}

Basic
equations of quantum theory are linear; this is the superposition
principle in quantum mechanics. The Hamilton--Jacobi equation,
the basic equation of classical mechanics, is nonlinear in the
conventional sense. However, it is linear over the semirings
$\rmax$ and $\rmin$. Similarly, different versions of the Bellman
equation, the basic equation of optimization theory, are linear
over suitable idempotent semirings; this is V.~P.~Maslov's idempotent
superposition principle, see \cite{Mas86,Mas87a,Mas87b,MaKo94,%
MaSa92}. For instance, the
finite-dimensional stationary Bellman equation can be written in
the form $X = H \odot X \oplus F$, where $X$, $H$, $F$ are
matrices with coefficients in an idempotent semiring $S$ and the
unknown matrix $X$ is determined by $H$ and $F$ \cite{Ca71}. In particular,
standard problems of dynamic programming and the well-known
shortest path problem correspond to the cases $S = \rmax$ and $S
=\rmin$, respectively. In \cite{Ca71}, it was known that principal 
optimization
algorithms for finite graphs correspond to standard methods for
solving systems of linear equations of this type (over
semirings). Specifically, Bellman's shortest path algorithm
corresponds to a version of Jacobi's algorithm, Ford's algorithm
corresponds to the Gauss--Seidel iterative scheme, etc.

The linearity of the Hamilton--Jacobi equation over
$\rmin$ and $\rmax$, which is the result of
the Maslov dequantization of the Schr{\"o}\-din\-ger equation, 
is closely related to the (conventional)
linearity of the Schr{\"o}\-din\-ger equation and can be deduced from
this linearity. Thus, it is possible to borrow standard ideas and
methods of linear analysis and apply them to a new area.

Consider a classical dynamical system specified by the Hamiltonian
$$
   H = H(p,x) = \sum_{i=1}^N \frac{p^2_i}{2m_i} + V(x),
$$
where $x = (x_1, \dots, x_N)$ are generalized coordinates, $p = (p_1,
\dots, p_N)$ are generalized momenta, $m_i$ are masses, and
$V(x)$ is the potential. In this case the Lagrangian $L(x, \dot x, t)$ has
the form
$$
   L(x, \dot x, t)
        = \sum^N_{i=1} m_i \frac{\dot x_i^2}2 - V(x),
$$
where $\dot x = (\dot x_1, \dots, \dot x_N)$, $\dot x_i = dx_i / dt$. The
value function $S(x,t)$ of the action functional has the form
$$
   S = \int^t_{t_0} L(x(t), \dot x(t), t)\, dt,
$$
where the integration is performed along the factual trajectory of the system.  The
classical equations of motion are derived as the stationarity conditions
for the action functional (the Hamilton principle, or the least action
principle).

For fixed values of $t$ and $t_0$ and arbitrary trajectories
$x(t)$,
the action functional $S=S(x(t))$ can be considered as a function
taking the set of curves (trajectories) to the set of real numbers
which can be treated as elements of  $\rmin$. In this case the
minimum of the action functional can be viewed as the Maslov
integral of this function over the set of trajectories or an
idempotent analog of the Euclidean version of the Feynman
path integral. The minimum of the
action functional corresponds to the maximum of $e^{-S}$, i.e.
idempotent integral $\int^{\oplus}_{\{paths\}} e^{-S(x(t))}
D\{x(t)\}$ with respect to the max-plus algebra $\rset_{\max}$.
Thus the least action principle can be considered as
an idempotent version of the well-known Feynman approach to
quantum mechanics.  The representation of a solution to the
Schr{\"o}\-din\-ger equation in terms of the Feynman integral 
corresponds to the Lax--Ole\u{\i}nik solution formula for the
Hamilton--Jacobi equation.

Since $\partial S/\partial x_i = p_i$, $\partial S/\partial t = -H(p,x)$,
the following Hamilton--Jacobi equation holds:
$$
   \pd{S}{t} + H \left(\pd{S}{x_i}, x_i\right)= 0.\eqno{(3)}
$$

Quantization (see, e.g., \cite{FeHi65}) leads to the Schr\"odinger equation
$$
   -\frac{\hbar}i \pd{\psi}{t}= \widehat H \psi = H(\hat p_i, \hat x_i)\psi,
   \eqno{(4)}
$$
where $\psi = \psi(x,t)$ is the wave function, i.e., a time-dependent
element of the Hilbert space $L^2(\rset^N)$, and $\widehat H$ is the energy
operator obtained by substitution of the momentum operators $\widehat p_i
= {\hbar \over i}{\partial \over \partial x_i}$ and the coordinate
operators $\widehat x_i \colon \psi \mapsto x_i\psi$ for the variables
$p_i$ and $x_i$ in the Hamiltonian function, respectively. This equation is
linear in the conventional sense (the quantum superposition principle). The
standard procedure of limit transition from the Schr\"odinger equation to
the Hamilton--Jacobi equation is to use the following ansatz for the wave
function:  $\psi(x,t) = a(x,t) e^{iS(x,t)/\hbar}$, and to keep only the
leading order as $\hbar \to 0$ (the `semiclassical' limit).

Instead of doing this, we switch to imaginary values of the Planck constant
$\hbar$ by the substitution $h = i\hbar$, assuming $h > 0$. Thus the
Schr\"odinger equation~(4) turns to the following generalized heat equation:
$$
   h\pd{u}{t} = H\left(-h\frac{\partial}{\partial x_i}, \hat x_i\right) u,
   \eqno{(5)}
$$
where the real-valued function $u$ corresponds to the wave function
$\psi$ (or, rather, to~$|\psi|$). A similar idea (the switch to imaginary time) is used in the
Euclidean quantum field theory (see, e.g., \cite{Ne73}); let us remember
that time and energy are dual quantities.

Linearity of equation~(4) implies linearity of equation~(5). Thus if
$u_1$ and $u_2$ are solutions of~(5), then so is their linear
combination
$$
   u = \lambda_1 u_1 + \lambda_2 u_2.\eqno{(6)}
$$

Let $S = h \ln u$ or $u = e^{S/h}$ as in Section 2 above. It can easily
be checked that equation~(5) thus turns to
$$
   \pd{S}{t}= V(x) + \sum^N_{i=1} \frac1{2m_i}\left(\pd{S}{x_i}\right)^2
        + h\sum^n_{i=1}\frac1{2m_i}\frac{\partial^2 S}{\partial x^2_i}.
   \eqno{(7)}
$$
Thus we have a passage from (4) to (7) by means of the change of variables
$\psi = e^{S/h}$.
Note that $|\psi| = e^{\Re S/h}$ , where $\Re S$ is the real part  of $S$. 
Now let us consider $S$ as a real variable. The equation 
(7) is nonlinear in the conventional sense. However, if $S_1$ and $S_2$
are its solutions, then so is the function
$$
   S = \lambda_1 \odot S_1 \opl_h \lambda_2\odot S_2
$$
obtained from~(6) by means of our substitution $S = h \ln u$.  Here the
generalized multiplication $\odot$ coincides with the ordinary addition and
the generalized addition $\opl_h$ is the image of the conventional addition
under the above change of variables.  As $h \to 0$, we obtain the
operations of the idempotent semiring $\rmax$, i.e., $\oplus = \max$ and
$\odot = +$, and equation~(7) turns to the Hamilton--Jacobi
equation~(3), since the third term in the right-hand side of
equation~(7) vanishes.

Thus it is natural to consider the limit function $S = \lambda_1 \odot S_1
\oplus \lambda_2 \odot S_2$ as a solution of the Hamilton--Jacobi equation and
to expect that this equation can be treated as linear over $\rmax$. This
argument (clearly, a heuristic one) can be extended to equations of a more
general form. For a rigorous treatment of (semiring) linearity for these
equations see \cite{FlMc2000,KoMa97,MaKo94,MaSa92} and also 
\cite{Mas87a}. Notice that if $h$ is changed
to $-h$, then we have that the resulting Hamilton--Jacobi equation is linear 
over $\rmin$.

The idempotent superposition principle indicates that there exist important
nonlinear (in the traditional sense) problems that are linear over idempotent 
semirings. The liner idempotent functional analysis (see below) is a
natural tool for investigation of those nonlinear 
infinite-dimensional problems that possess this property.

\section{Convolution and the Fourier--Legendre transform}

Let $G$ be a group. Then the space $\CalB(X, \rset_{\max})$ of all bounded
functions $G\to\rset_{\max}$ (see above) is an idempotent semiring
with respect to the following analog $\circledast$ of the usual 
convolution:
$$
   (\varphi(x)\circledast\psi)(g)=
         \int_G^{\oplus} \varphi (x)\odot\psi(x^{-1}\cdot g)\, dx=
\sup_{x\in G}(\varphi(x)+\psi(x^{-1}\cdot g)).
$$
Of course, it is possible to consider other ``function spaces'' (and other
basic semirings instead of $\rset_{\max}$). In \cite{KoMa97,MaKo94} ``group
semirings'' of this type are referred to as {\it convolution
semirings}.

Let $G=\rset^n$, where $\rset^n$ is considered as a topological group
with respect to the vector addition. The conventional Fourier--Laplace 
transform is defined as
$$
   \varphi(x) \mapsto \tilde{\varphi}(\xi)
        = \int_G e^{i\xi \cdot x} \varphi (x)\, dx,\eqno{(8)}
$$
where $e^{i\xi \cdot x}$ is a character of the group $G$, i.e., a solution
of the following functional equation:
$$
   f(x + y) = f(x)f(y).
$$
The idempotent analog of this equation is
$$
   f(x + y) = f(x) \odot f(y) = f(x) + f(y),
$$
so ``continuous idempotent characters'' are linear functions of the
form $x \mapsto \xi \cdot x = \xi_1 x_1 + \dots + \xi_n x_n$. As a result,
the transform in~(8) assumes the form
$$
   \varphi(x) \mapsto \tilde{\varphi}(\xi)
        = \int_G^\oplus \xi \cdot x \odot \varphi (x)\, dx
   = \sup_{x\in G} (\xi \cdot x + \varphi (x)).\eqno{(9)}
$$
The transform in~(9) is nothing but the {\it Legendre transform\/}
(up to some notation) \cite{Mas87b}; transforms of this kind establish the
correspondence between the Lagrangian and the Hamiltonian formulations
of classical mechanics. The Legendre transform generates an idempotent
version of harmonic analysis for the space of convex functions,
see, e.g., \cite{MaTi2003}.

Of course, this construction can be generalized to different classes
of groups and semirings. Transformations of this type convert the
generalized convolution $\circledast$ to the pointwise (generalized)
multiplication and possess analogs of some important properties
of the usual Fourier transform. For the case of semirings of Pareto
sets, the corresponding version of the Fourier transform reduces the
multicriterial optimization problem to a family of singlecriterial 
problems \cite{SaTa92}.

The examples discussed in this sections can be treated as fragments
of an idempotent version of the representation theory, see, e.g.,
\cite{LiMaSh2002}. In particular,
``idempotent'' representations of groups can be examined as representations
of the corresponding convolution semirings (i.e. idempotent group semirings)
in semimodules.

%% {\bf 8. Correspondence to stochastics and a duality between 
%% probability and optimization.}
\section{Correspondence to stochastics and a duality between 
probability and optimization}

Maslov measures are nonnegative (in the sense of the standard
order) just as probability measures. The analogy between idempotent and
probability measures leads to important relations between
optimization theory and probability theory. By the present time
idempotent analogues of many objects of
stochastic calculus have been constructed and investigated, 
such as  max-plus martingales, max-plus
stochastic differential equations, and others. These results allow, for
example, to transfer powerful stochastic methods
to the optimization theory. This was noticed and examined by
many authors (G.~Salut, P.~Del Moral, M.~Akian, J.-P.~Quadrat,
V.~P.~Maslov, V.~N.~Kolokoltsov, P.~Bernhard, W.~A.~Fleming,
W.~M.~McEneaney, A.~A.~Puhalskii and others), see the survey paper of
W.~A.~Fleming and W.~M.~McEneaney \cite{FlMc2005} 
and \cite{Ak99,AkQuVi98,Be94,CoQu94,De97,De2004,DeDo98,%
DeDo2001,%
Fl2002,Fl2004,FlMc2000,Gun98a,MaKo94,Puh2001,Qu90,Qu94}.
For relations and applications to large deviations see \cite{Ak99,%
De97,De2004,DeDo98,DeDo2001,Puh2003} and
especially the book of A.~A.~Puhalskii \cite{Puh2001}.

%% {\bf 9. Idempotent functional analysis}. 
\section{Idempotent functional analysis}

Many other idempotent analogs may be given, in particular, for
basic constructions and theorems of functional analysis.
Idempotent functional
analysis is an abstract version of idempotent analysis. For the
sake of simplicity take $S=\rmax$ and let $X$ be an arbitrary set.
The idempotent integration can be defined by the formula (1), see
above. The functional $I(\varphi)$ is linear over $S$ and its
values correspond to limiting values of the corresponding analogs
of Lebesgue (or Riemann) sums. An idempotent scalar product of
functions $\varphi$ and $\psi$ is defined by the formula
$$
\langle\varphi,\psi\rangle = \int^{\oplus}_X \varphi(x)\odot\psi(x)\, dx =
\sup_{x\in X}(\varphi(x)\odot\psi(x)).
$$
So it is natural to construct idempotent analogs of integral
operators in the form
$$
\varphi(y) \mapsto (K\varphi)(x) = \int^{\oplus}_Y K(x,y)\odot
\varphi(y)\, dy = \sup_{y\in Y}\{K(x,y)+\varphi(y)\},\eqno(10)
$$
where $\varphi(y)$ is an element of a space of functions defined
on a set $Y$, and $K(x,y)$ is an $S$-valued function on $X\times Y$.
Of course, expressions of this type are standard in optimization
problems.\medskip

Recall that the definitions and constructions described above can
be extended to the case of idempotent semirings which are
conditionally complete in the sense of the standard order. Using
the Maslov integration, one can construct various function spaces
as well as idempotent versions of the theory of generalized
functions (distributions). For some concrete idempotent function
spaces it was proved that every `good' linear operator (in the
idempotent sense) can be presented in the form (10); this is an
idempotent version of the kernel theorem of L.~Schwartz; results
of this type were proved by V.~N.~Kolokoltsov, P.~S.~Dudnikov and
S.~N.~Samborski\u\i, I.~Singer, M.~A.~Shubin and others, see, e.g., 
\cite{DuSa92,KoMa97,MaKo94,MaSa92,Si2003}. So every
`good' linear functional can be presented in the form
$\varphi\mapsto\langle\varphi,\psi\rangle$, where
$\langle,\rangle$ is an idempotent scalar product.\medskip

In the framework of idempotent functional analysis results of this
type can be proved in a very general situation. In \cite{LiMaSh98,%
LiMaSh99,LiMaSh2001,LiMaSh2002,LiSh2002} an algebraic version 
of the idempotent functional
analysis is developed; this means that basic (topological) notions
and results are simulated in purely algebraic terms. The treatment
covers the subject  from basic concepts and results (e.g.,
idempotent analogs of the well-known theorems of Hahn-Banach,
Riesz, and Riesz-Fisher) to idempotent analogs of 
A.~Grothendieck's concepts and results on topological tensor
products, nuclear spaces and operators. An abstract version of the
kernel theorem is formulated. Note that the passage from the usual
theory to idempotent functional analysis may be very nontrivial;
for example, there are many non-isomorphic idempotent Hilbert
spaces. Important results on idempotent functional analysis
(duality and separation theorems) are recently published by 
G.~Cohen, S.~Gaubert, and J.-P.~Quadrat \cite{CoGaQu2004};
see also a finite dimensional version of the separation
theorem in \cite{ZiK77}.
Idempotent functional analysis has received much attention
in the last years, see, e.g., \cite{AkGa2003,AkGaKo2005,AkGaWa2004,%
AkGaWa2005,CoGaQuSi2005,GoMi2001,LiMa2005,LiSh2005,%
Rou2005,Si2003,So99a,Wal2005} and works cited above.

%%%\section{10. The dequantization transform and generalization of the
%%%Newton polytopes}
\section{The dequantization transform and generalization of the
Newton polytopes}

In this section we briefly discuss results proved in the paper
\cite{LiSh2005}. For functions defined on $\cset^n$ almost everywhere (a.w.)
it is possible to construct a dequantization transform $f\to\hat f$
generated by the Maslov dequantization. If $f$ is a polynomial, then the
subdifferential $\partial{\hat f}$ of $\hat f$ at the origin coincides with
the Newton polytope of $f$. For the semiring of polynomials with nonnegative 
coefficients, the dequantization transform is a homomorphism of this semiring
to the idempotent semiring of convex polytopes with respect to the
well-known Minkovski operations. These results can be generalized to
a wide class of functions and convex sets.

\subsection{The dequntization transform}
Let $X$ be a topological space. For functions $f(x)$ defined on $X$
we shall say that a certain assertion is valid {\it almost everywhere} (a.e.) if 
it is valid for all elements $x$ of an open dense subset of $X$. 
Suppose $X$ is $\cset^n$ or $\rset^n$; denote by $\rset^n_+$ the set
$x=\{\,(x_1, \dots, x_n)\in X \mid x_i\geq 0$ for 
$i = 1, 2, \dots, n\}$. For $x= (x_1, \dots, x_n)\in X$ we set 
$\exp(x) = (\exp(x_1), \dots, \exp(x_n))$;
so if $x\in\rset^n$, then $\exp(x)\in \rset^n_+$. 

Denote by $\maF(\cset^n)$ the set of all functions defined
and continuous on an open dense subset $U\subset \cset^n$
such that $U\supset \rset^n_+$. In all the examples below we consider
even more regular functions, which are holomorphic in~$U$.
 It is clear that $\maF(\cset^n)$ is a ring 
(and an algebra over $\cset$) with respect to the usual addition
and multiplications of functions.

For $f\in \maF(\cset^n)$ let us define the function $\hat f_h$
by the following formula:
$$
\hat f_h(x) = h \log|f(\exp(x/h))|,\eqno(11)
$$
where $h$ is a (small) real parameter and $x\in\rset^n$. Set
$$
\hat f(x) = \lim_{h\to 0} \hat f_h (x),\eqno(12)
$$
if the right-hand part of (12) exists almost everywhere. We shall say that the 
function $\hat f(x)$ is a {\it dequantization} of the function $f(x)$ and the 
map $f(x)\mapsto \hat f(x)$ is a {\it dequantization transform}. By 
construction, $\hat f_h(x)$ and $\hat f(x)$ can be treated as functions taking their 
values in $\rset_{\max}$. Note that in fact $\hat f_h(x)$ and $\hat f(x)$
 depend on the restriction of $f$ to $\rset_+^n $
only; so in fact the dequantization transform is constructed
for functions defined on $\rset^n_+$ only. 
It is clear that the dequantization transform is 
generated by the Maslov dequantization and the map $x\mapsto |x|$. Of 
course, similar definitions can be given for functions defined on $\rset^n$
and $\rset_+^n$.

Denote by $V$ the set $\rset^n$ treated as a linear Euclidean
space (with the scalar product $(x, y) = x_1y_1+ x_2y_2 +\dots + x_ny_n$)
and set $V_+ = \rset_+^n$.
We shall say that a function $f\in \maF(\cset^n)$ is {\it dequantizable} 
whenever its 
dequantization $\hat f(x)$ exists (and is defined on an open dense subset 
of $V$). By $\maD (\cset^n)$ denote the set of all dequantizable functions and by 
$\widehat{\maD}(V)$ denote the set $\{\,\hat f \mid f\in \maD(\cset^n)\,\}$.
 Recall that 
functions from $\maD(\cset^n)$ (and $\widehat{\maD}(V)$) are defined almost
 everywhere and 
$f=g$ means that $f(x) = g(x)$ a.e., i.e., for $x$ ranging over an open dense subset 
of $\cset^n$ (resp., of $V$). Denote by $\maD_+(\cset^n)$ the set of all 
functions $f\in \maD(\cset^n)$ 
such that $f(x_1, \dots, x_n)\geq 0$ if $x_i\geq 0$ for $i= 1,\dots, n$; so
 $f\in \maD_+(\cset^n)$ if the restriction of $f$ to $V_+ = \rset_+^n$ is a
 nonnegative function. By $\widehat{\maD}_+(V)$ denote the image of $\maD_+(\cset^n)$
 under the dequantization transform. We shall say that functions 
$f, g\in \maD(\cset^n)$ are in {\it general position} whenever
$\hat f (x) \neq \widehat g(x)$ for $x$ running an open dense
subset of $V$.

For functions $f, g \in \maD(\cset^n)$ and any nonzero constant $c$, the 
following equations are valid:

\begin{enumerate}
\item[1)] $\widehat{fg} = \hat f + \widehat g$;
\item[2)] $|\hat f| = \hat f$; $\widehat{cf} = f$; $\widehat c =0$;
\item[3)] $(\widehat{f+g})(x) = \max\{\hat f(x), \widehat g(x)\}$ a.e.\ if $f$ and $g$
are nonnegative on $V_+$ (i.e., $f, g \in \maD_+(\cset^n)$) or $f$ and $g$
 are in general position.
\end{enumerate}
Left-hand sides of these equations are well-defined automatically.

The set $\maD_+(\cset^n)$ has a natural structure of a semiring with respect to 
the
 usual addition and multiplication of functions taking their values in $\cset$.
 The set $\widehat{\maD}_+(V)$ has a natural
 structure of an idempotent semiring with respect to the operations 
$(f\oplus g)(x) = \max \{ f(x), g(x)\}$, $(f\odot g)(x) = f(x) + g(x)$; 
elements
 of $\widehat{\maD}_+(V)$ can be naturally treated as functions taking their values
 in $\rset_{\max}$. The dequantization transform generates a homomorphism from
 $\maD_+(\cset^n)$ to $\widehat{\maD}_+(V)$.

\subsection{Simple functions}
For any nonzero number $a\in\cset$ and any vector
 $d = (d_1, \dots, d_n)\in V = \rset^n$ 
we set $m_{a,d}(x) = a \prod_{i=1}^n x_i^{d_i}$; functions of this kind we
 shall 
call {\it generalized monomials}. Generalized monomials are defined a.e.\ on
$\cset^n$ and on $V_+$, but not on $V$ unless the numbers $d_i$ take
integer or suitable rational values. We shall say that a function $f$ is a {\it 
generalized polynomial} whenever it is a finite sum of linearly
 independent generalized monomials. For instance, Laurent polynomials
are examples of generalized polynomials.

As usual, for $x, y\in V$ we set $(x,y) = x_1y_1 + \dots + x_ny_n$. It is 
easy to prove that if $f$ is a generalized monomial 
$m_{a,d}(x)$, then $\hat f$ is a linear function $x\mapsto (d,x)$.
If $f$ is a generalized polynomial, then $\hat f$ is a sublinear function.

Recall that a real function $p$ defined on $V = \rset^n$ is {\it sublinear} 
if $p = \sup_{\alpha} 
p_{\alpha}$, where $\{p_{\alpha}\}$ is a collection of linear functions. 
Sublinear functions defined everywhere on $V=\rset^n$ are convex; thus these
 functions are continuous.
We discuss sublinear functions of this kind only. Suppose $p$ is a continuous
 function defined on $V$, then $p$ is sublinear whenever 

1) $p(x+ y) \leq p(x) + p(y)$ for all $x, y \in V$;

2) $p(cx) = cp(x)$ for all $x\in V$, $c\in \rset_+$.

So if $p_1$, $p_2$ are sublinear functions, then $p_1 +p_2$ is a sublinear
 function.

We shall say that a function $f \in \maF(\cset^n)$ is {\it simple}, if 
its dequantization $\hat f$ exists and a.e.\ coincides with a sublinear 
function; by misuse of language, we shall denote this (uniquely defined verywhere on $V$) sublinear function by the same symbol $\hat f$.

Recall that simple functions $f$ and $g$ are {\it in general position} if 
$\hat f(x) \neq \widehat g(x)$ for all $x$ belonging to an open dense subset of $V$. In particular, generalized monomials are in  
general position whenever they are linearly independent.

Denote by $\mathit{Sim}(\cset^n)$ the set of all simple functions defined on 
$V$ and denote by $\mathit{Sim}_+(\cset^n)$ the set $\mathit{Sim}(\cset^n) 
\cap \maD_+(\cset^n)$. By $\mathit{Sbl}(V)$ denote the
 set of all (continuous) sublinear functions defined on $V = \rset^n$ and by
 $\mathit{Sbl}_+(V)$ denote the image $\widehat{\mathit{Sim}_+}(\cset^n)$ of 
$\mathit{Sim}_+(\cset^n)$ under the dequantization transform.

The set $\mathit{Sim}_+(\cset^n)$ is a subsemiring of $\maD_+(\cset^n)$ and $\mathit{Sbl}_+(V)$
is an idempotent subsemiring of $\widehat{\maD_+}(V)$. The
dequantization transform generates an epimorphism of
$\mathit{Sim}_+(\cset^n)$ onto $\mathit{Sbl}_+(V)$. The set $\mathit{Sbl}(V)$ is an idempotent
semiring with respect to the operations
$(f\oplus g)(x) = \max \{ f(x), g(x)\}$, 
$(f\odot g)(x) = f(x) + g(x)$.

Of course, polynomials and generalized polynomials are simple functions.

We shall say that functions $f, g\in\maD(V)$ are {\it asymptotically equivalent} 
whenever $\hat f = \widehat g$; any simple function $f$ is an {\it asymptotic 
monomial} whenever 
$\hat f$ is a linear function. A simple function $f$ will be called an {\it 
asymptotic polynomial} whenever $\hat f$ is a sum of a finite collection of 
nonequivalent asymptotic monomials.
Every asymptotic polynomial is a simple function.

{\bf Example.} Generalized polynomials, logarithmic functions of 
(generalized) polynomials, and products of 
polynomials and logarithmic functions are asymptotic polynomials. This follows 
from our definitions and formula~(11).

\subsection{Subdifferentials of sublinear functions and Newton sets
for simple functions}
It is well known that all the convex compact subsets in $\rset^n$ form an 
idempotent 
semiring $\mathcal{S}$ with respect to the Minkowski operations: for $A, B \in \mathcal{S}$ 
the sum $A\oplus B$ is the convex hull of the union $A\cup B$; the product 
$A\odot B$ is defined in the following way: $A\odot B = \{\, x\mid x = a+b, 
\text{ where $a\in A, b\in B$}\,\}$. In fact $\mathcal{S}$ is an idempotent
 linear space over $\rset_{\max}$ (see, e.g., \cite{LiMaSh2001}). Of course, 
the Newton polytopes in $V$ form a 
subsemiring $\mathcal{N}$ in $\mathcal{S}$. 

We shall use some elementary results from convex analysis. These results can be 
found, e.g., in \cite{MaTi2003}.

For any function $p\in \mathit{Sbl}(V)$ we set
$$
\partial p = \{\, v\in V\mid (v, x) \le p(x)\ \forall x\in V\,\}.
$$

It is well known from convex analysis that for any sublinear function $p$ the 
set $\partial p$ is exactly the {\it subdifferential} of $p$ at the origin. 
 The following proposition is also known in convex 
analysis.

\begin{prop}
Suppose $p_1,p_2\in \mathit{Sbl}(V)$, then
\begin{center}
$\partial (p_1+p_2) = \partial p_1\odot\partial p_2 =$

$ \{\, v\in V\mid
v = v_1+v_2, \text{ where $v_1\in \partial p_1, v_2\in \partial p_2$}\,\}$;

$\partial (\max\{p_1(x), p_2(x)\}) = \partial p_1\oplus\partial p_2$.
\end{center}

Suppose $p\in \mathit{Sbl}(V)$.  Then $\partial p$ is a nonempty convex compact 
subset of $V$.
\end{prop}

\begin{cor}
The map $p\mapsto \partial p$ is a homomorphism of the idempotent semiring
 $\mathit{Sbl}(V)$ to the idempotent semiring $\mathcal{S}$ of all convex
 compact subsets of $V$.
\end{cor}

For any simple function $f\in \mathit{Sim}(\cset^n)$ let us denote by $N(f)$ 
the set $\partial(\hat f)$. We shall call $N(f)$ the {\it Newton set} of the
 function $f$. It follows from this proposition that for 
any simple function $f$, its Newton set $N(f)$ is a nonempty convex
 compact subset of $V$.

\begin{thm}
Suppose that $f$ and $g$ are simple functions. Then
\begin{enumerate}
\item[1)] $N(fg) = N(f)\odot N(g) = \{\, v\in V\mid v = v_1 +v_2$
with  $v_1 \in N(f), v_2 \in N(g)\}$;
\item[2)] $N(f+g) = N(f)\oplus N(g)$, if $f_1$ and $f_2$ are in general
 position or $f_1, f_2 \in \mathit{Sim}_+(\cset^n)$ {\rm (}recall that $N(f)\oplus N(g)$ 
is the convex hull of $N(f)\cup N(g)${\rm )}.
\end{enumerate}
\end{thm}

\begin{cor}
The map $f\mapsto N(f)$ generates a homomorphism from $\mathit{Sim}_+(\cset^n)$ to 
$\mathcal{S}$.
\end{cor}

\begin{prop}
Let $f = m_{a,d}(x) = a \prod^n_{i=1} x_i^{d_i}$ be a generalized monomial; here
$d = (d_1, \dots, d_n) \in V= \rset^n$ and $a$ is a nonzero
complex number. Then $N(f) = \{ d\}$.
\end{prop}
 
\begin{cor}
Let $f = \sum_{d\in D} m_{a_d,d}$ be a generalized polynomial. Then $N(f)$ is the polytope 
$\oplus_{d\in D}\{d\}$, i.e.\ the convex hull of the finite set $D$.
\end{cor}

This statement follows from Theorem and Proposition 2. Thus in this case
 $N(f)$ is the well-known classical Newton polytope of the polynomial $f$.

Now the following corollary is obvious.

\begin{cor}
Let $f$ be a generalized or asymptotic polynomial. Then its Newton set
 $N(f)$ is a convex polytope.
\end{cor}

{\bf Example.} Consider the one dimensional case, i.e., $V = \rset$ and 
suppose
 $f_1 = a_nx^n + a_{n-1}x^{n-1} + \dots + a_0$ and $f_2 = b_mx^m + b_{m-1}
 x^{m-1} + \dots + b_0$, where $a_n\neq 0$, $b_m\neq 0$, $a_0 \neq 0$,
 $b_0 \neq 0$. Then $N(f_1)$ is the segment $[0, n]$ and $N(f_2)$ is the
 segment $[0, m]$. So the map $f\mapsto N(f)$ corresponds to the map
 $f\mapsto \deg (f)$, where $\deg(f)$ is a degree of the polynomial $f$. In
 this case Theorem means that $\deg(fg) = \deg f + \deg g$ and
 $\deg (f+g) = \max \{\deg f, \deg g\} = \max \{n, m\}$ if $a_i\geq 0$,
 $b_i\geq 0$ or $f$ and $g$ are in general position.

%%{\bf 11. Dequantization of geometry.} 
\section{Dequantization of geometry}

An idempotent version of real
algebraic geometry was discovered in the report of O.~Viro for the
Barcelona Congress \cite{Vir2000}. Starting from the idempotent 
correspondence principle O.~Viro constructed a piecewise-linear
geometry of polyhedra of a special kind in finite dimensional 
Euclidean spaces as a result of the
Maslov dequantization of real algebraic geometry. He indicated
important applications in real algebraic geometry (e.g., 
 in the framework of Hilbert's 16th problem for
constructing real algebraic varieties with prescribed properties
and parameters) and relations to complex algebraic geometry and
amoebas in the sense of I.~M.~Gelfand, M.~M.~Kapranov, and
A.~V.~Zelevinsky (see their book \cite{GeKaZe94} and \cite{Vir2002}).
Then complex algebraic geometry was dequantized by 
G.~Mikhalkin and the result turned out
to be the same; this new `idempotent' (or
asymptotic) geometry is now often called the {\it tropical algebraic geometry},
see, e.g., \cite{EiKaLi2004,ItKhSh2003,Mi2001,Mi2003,Mi2004,%
Mi2005,RiStTh2005,Shu2002,SpSt2004,St2002}.

There is a natural relation between the Maslov dequantization
and amoebas. 
Suppose $({\cset}^*)^n$ is a complex torus, where
${\cset}^* = {\cset}\backslash \{0\}$ is the group
of nonzero complex numbers under multiplication.  For
 $z = (z_1, \dots, z_n)\in 
(\cset^*)^n$ and a positive real number $h$ denote by
$\Log_h(z) = h\log(|z|)$ the element
\[(h\log |z_1|, h\log |z_2|, \dots,
h\log|z_n|) \in \rset^n.\]  
Suppose $V\subset (\cset^*)^n$ is a
complex algebraic variety; denote by $\maA_h(V)$ the set
$\Log_h(V)$. If $h=1$, then the set $\maA(V) = \maA_1(V)$ is
called the {\it amoeba} of $V$ in the sense of \cite{GeKaZe94},
see also \cite{Ale2005,Mi2001,Mi2004,Mi2005,PaTs2005,%
Shu2002,So2005,Vir2002};
the amoeba $\maA(V)$ is a closed subset of $\rset^n$ with a
non-empty complement. 
Note that this construction depends on our coordinate
system.

For the sake of simplicity suppose $V$ is a hypersurface 
in~$(\cset^*)^n$ defined by a polynomial~$f$; then there 
is a deformation $h\mapsto f_h$ of this polynomial generated by the
Maslov dequantization and $f_h = f$ for $h = 1$.
Let $V_h\subset ({\cset}^*)^n$ be the zero set of $f_h$ and
set $\maA_h (V_h) = {\Log}_h (V_h)$. Then
 there exists a tropical variety
$\mathit{Tro}(V)$ such that the subsets $\maA_h(V_h)\subset \rset^n$
tend to $\mathit{Tro}(V)$ in the Hausdorff metric as $h\to 0$, 
see \cite{Mi2001,Ru2001}.
The tropical variety $\mathit{Tro}(V)$ is a result
of a deformation of the amoeba $\maA(V)$ and the Maslov
dequantization of the variety $V$. The set $\mathit{Tro}(V)$ is called the
{\it skeleton} of $\maA(V)$. 

\begin{figure}[t]
  \label{fig:amoeba}
  \centering
  \includegraphics{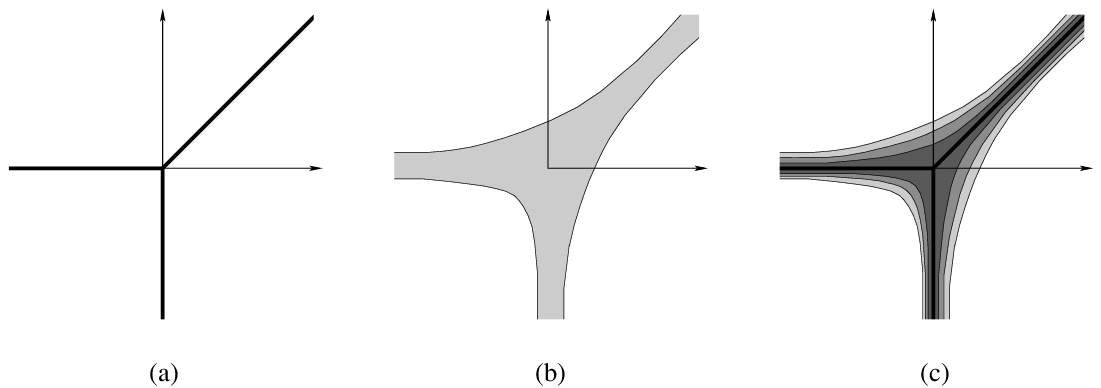}
  \caption{A tropical line and amoebas.}
\end{figure}

{\bf Example \cite{Mi2001}.}  For the line $V = \{\, (x, y)\in
({\cset}^*)^2 \mid x + y + 1 = 0\,\}$ the piecewise-linear graph
$\mathit{Tro}(V)$ is a tropical line, see Fig.~\ref{fig:amoeba}(a).
The amoeba $\maA(V)$ is represented in Fig.~\ref{fig:amoeba}(b), while
Fig.~\ref{fig:amoeba}(c) demonstrates the corresponding deformation of
the amoeba.

In the important paper \cite{Kap2000} (see also 
\cite{EiKaLi2004,Mi2001,Mi2004,RiStTh2005}) tropical
varieties appeared as amoebas over non-Archimedean fields. In
2000 M.~Kontsevich noted that it is possible to use non-Arhimedean
amoebas in enumerative geometry, see \cite[section
2.4, remark 4]{Mi2001}. In fact methods of tropical geometry lead to
remarkable applications to the algebraic enumerative geometry,
Gromov-Witten and Welschinger invariants, see \cite{GaMa2005a,ItKhSh2003,%
ItKhSh2003a,ItKhSh2004,ItKhSh2005,Mi2001,Mi2003,Mi2004,Mi2005,Shu2002,Shu2004}.
In particular,
G.~Mikhalkin presented and proved in \cite{Mi2003,Mi2005} a
formula enumerating curves of arbitrary genus in toric
surfaces. See also the papers \cite{GaMa2005b,ItKhSh2003,ItKhSh2003a,NiSie2004,Shu2002}.

Last time many other papers on tropical algebraic geometry and its applications
to the conventional (e.g., complex) algebraic geometry and other areas appeared,
see, e.g., \cite{Ale2005,FoGo2003a,FoGo2003b,Izh2005a,NiSie2004,PaSt2003,PaSt2004,%
SzVe2004,Te2004,Vi2004}. The thing is that some difficult traditional problems
can be reduced to their tropical versions which are hopefully not so difficult.

Note that tropical geometry is closely related to the
well-known program of M.~Kontsevich and Y.~Soibelman, see,
e.g., \cite{KonSoi2001,KonSoi2004}. 

There is an introductory
paper \cite{RiStTh2005} (see also \cite{St2002}) on tropical algebraic geometry.
However, on the whole, only first steps in idempotent/tropical geometry
have been made and the problem of systematic construction of idempotent
versions of algebraic and analytic geometries is still open.

%% {\bf 12. Correspondence principle for algorithms and their computer
%% implementations.} 
\section{The correspondence principle for algorithms and their computer
implementations}

There are many important applied algorithms of idempotent
mathematics, see, e.g., \cite{BaCoOlQu92,Bu2005,Ca79,Cu95,CuMe80,FiRo93,%
FlMc2000,FlMc2005,GoMi79,GoMi2001,Gun98a,ItKhSh2003,KiRo2004,Ki2001,%
Ko2001,KoMa97,LiMa95,LiMa98,LiMaE2000,LiMaRo2000,LiSo2000,LiSo2001,%
LoPe2005,Mi2003,Mi2005,RiStTh2005,Roy92,St2002,Vor63,Vor67,Vor70,% 
WoOl2005,ZiK2005,ZiU81}. 

The~idempotent correspondence principle is valid
for algorithms as well as for their software and hardware implementations
\cite{LiMa95,LiMa96,LiMa98,LiMaE2000,LiMaRo2000}. 
In particular, according to the superposition principle, analogs
of linear algebra algorithms are especially important. It is well-known
that algorithms of linear algebra are convenient for parallel computations;
so their idempotent analogs accept a parallelization. This is a regular way
to use parallel computations for many problems including basic optimization
problems. It is convenient to use universal algorithms which do not
depend on a concrete semiring and its concrete computer model. Software
implementations for universal semiring algorithms are based on
object-oriented and generic programming; program modules can deal with
abstract (and variable) operations and data types,
see \cite{LiMa95,LiMa98,LiMaE2000,LiMaRo2000,LoPe2005}. 

The most important and standard algorithms have many hardware
implementations in the form
of technical devices or special processors. These devices often can be used as
prototypes for new hardware units generated by substitution of the
usual arithmetic operations for its semiring analogs, see 
\cite{LiMa95,LiMa98,LiMaRo2000}.
Good and efficient technical ideas and decisions can be
transposed from prototypes into new hardware units. Thus the correspondence
principle generates a regular heuristic method for hardware design.

%%{\bf 13. Idempotent interval analysis.}
\section{Idempotent interval analysis}

An idempotent version of the
traditional interval analysis is presented in \cite{LiSo2000,LiSo2001}.
Let $S$ be an idempotent semiring equipped with the
standard partial order. A {\it closed interval} in $S$ is a subset
of the form ${\bf x} = [{\bf \underline x}, {\bf \bar  x}] = \{
x\in S\mid {\bf \underline x}\cle x
\cle{\bf \bar x}\}$, where the elements $\bf \underline x\cle \bf
\bar x$ are called {\it lower} and {\it upper bounds} of the
interval $\bf x$. A {\it weak interval extension} $I(S)$ of the
semiring $S$ is the set of all closed intervals in $S$ endowed
with operations $\oplus$ and $\odot$ defined as $\bf x\oplus\bf y
= [\bf \underline x \oplus\bf \underline y, \bf \bar
x\oplus\bf\bar y]$, $\bf x\odot\bf y = [\bf \underline x\odot\bf
\underline y, \bf \bar x\odot\bf \bar y]$; the set $I(S)$ is a new
idempotent semiring with respect to these operations.  It is
proved that basic interval problems of idempotent linear algebra
are polynomial, whereas in the traditional interval analysis
problems of this kind are generally NP-hard. Exact interval
solutions for the discrete stationary Bellman equation (see the
matrix equation discussed in section 8 above) and for the 
corresponding optimization problems are constructed and examined by
G.~L.~Litvinov and A.~N.~Sobolevski{\u\i}
in \cite{LiSo2000,LiSo2001}. Similar results are presented by
K.~Cechl{\'a}rov{\'a} and R.~A.~Cu\-ning\-ha\-me-Gre\-en
in~\cite{CeCu2002}.

%%{\bf 14. Relations to the KAM theory and optimal transport.} 
\section{Relations to the KAM theory and optimal transport}

The subject of the Kolmogorov--Arnold--Moser (KAM) theory may be
formulated as the study of invariant subsets in phase spaces of
nonintegrable Hamiltonian dynamical systems where the dynamics
displays the same degree of regularity as that of integrable systems
(quasiperiodic behaviour). Recently, a considerable progress was made
via a variational approach, where the dynamics is specified by the
Lagrangian rather than Hamiltonian function. The corresponding theory
was initiated by S.~Aubry and J.~N.~Mather and recently dubbed {\it
  weak KAM theory} by A.~Fathi (see his monograph ``Weak KAM Theorems
in Lagrangian Dynamics,'' in preparation, and also
\cite{KKS2005a,KKS2005b,So99a,So99b}).  Minimization of a certain
functional along trajectories of moving particles is a central feature
of another subject, the optimal transport theory, which also has
undergone a rapid recent development. This theory dates back to
G.~Monge's work on cuts and fills (1781). A modern version of the
theory is known now as the so-called {\it Monge--Amp\`ere--Kantorovich
  (MAK) optimal transport theory} (after the work of L.~Kantorovich
\cite{K42}).  There is a similarity between the two theories and there
are relations to problems of the idempotent functional analysis (e.g.,
the problem of eigenfunctions for ``idempotent'' integral operators,
see \cite{So99a}). Applications of optimal transport to data
processing in cosmology are presented in \cite{BFHLMMS2003,FMMS2002}.

%%{\bf 15. Relations to logic, fuzzy sets, and possibility theory.}
\section{Relations to logic, fuzzy sets, and possibility theory}

Let $S$ be an idempotent semiring with neutral elements $\0$ and
$\1$ (recall that $\0\neq \1$, see section 2 above). Then the Boolean
algebra $\bset = \{\0, \1\}$ is a natural idempotent subsemiring
of $S$. Thus $S$ can be treated as a generalized (extended) logic
with logical operations $\oplus$ (disjuction) and $\odot$ 
(conjunction). Ideas of this kind are discussed in many books and
papers with respect to generalized versions of logic and especially
 quantum logic, see, e.g., \cite{DiGe2005,Gol99,KlPa2004,Ro90,Ro96}.

Let $\Omega$ be the so-called universe consisting of ``elementary events.''
Denote by ${\CF}(S)$ the set of functions defined on $\Omega$ and taking
their values in $S$; then ${\CF}(S)$ is an idempotent semiring with respect
to the pointwise addition and multiplication of functions. We shall
say that elements of ${\CF}(S)$ are {\it generalized fuzzy sets}, 
see \cite{Gol99,Li2004}. We have the 
well-known classical definition of fuzzy sets (L.~A.~Zadeh \cite{Za65})
if $S = {\pset}$, where $\pset$ is the segment $[0,1]$ with the semiring operations
$\oplus = \max$ and $\odot = \min$. Of course,
functions from ${\CF}(\pset)$ taking their values in the Boolean algebra
$\bset = \{0, 1\}\subset {\pset}$ correspond to traditional sets from $\Omega$
and semiring operations correspond to standard operations for sets. In the
general case 
functions from ${\CF}(S)$ taking their values in $\bset = \{{\0}, 
{\1}\}\subset S$ can be treated as traditional subsets in $\Omega$. If $S$
is a lattice (i.e. $x\odot y = \inf \{x, y\}$ and $x\oplus y = \sup 
\{x, y\}$), then generalized fuzzy sets coincide with $L$-fuzzy sets in the
sense of J.~A.~Goguen \cite{Go67}. The set $I(S)$ of intervals is an idempotent
semiring (see section 11), so elements of ${\CF}(I(S))$ can be treated as
interval (generalized) fuzzy~sets.

It is well known that the classical theory of fuzzy sets is a basis
for the theory of possibility \cite{Za78,DuPrSa2001}. Of course, it is possible
to develop a similar generalized theory of possibility starting from
generalized fuzzy sets, see, e.g., \cite{DuPrSa2001,KlPa2004,Li2004}.
Generalized theories can be 
noncommutative; they seem to be more qualitative and less quantitative
with respect to the classical theories presented in \cite{Za65,Za78}. 
We see that
idempotent analysis and the theories of (generalized) fuzzy sets and 
possibility have the same objects, i.e. functions taking their values in
semirings. However, basic problems and methods could be different
for these theories (like for the measure theory and the probability
theory).

%%{\bf 16. Relations to other areas and miscellaneous applications.}
\section{Relations to other areas and miscellaneous applications}

Many relations and applications of idempotent mathematics to different
theoretical and applied areas of mathematical sciences are discussed
above. Of course, optimization and optimal control problems form a very natural 
field for applications of ideas and methods of idempotent mathematics. There
is a very good survey paper \cite{Ko2001} of V.~N.~Kolokoltsov on the
subject, see also \cite{BaCoOlQu92,Bu2005,Ca79,CoGaQu99,CoQu94,Cu79,Cu95,%
CuMe80,De97,DeDo98,DeDo2001,Fl2002,Fl2004,FlMc2000,FlMc2005,%
GoMi79,GoMi2001,Gun98a,LiMa95,LiMa98,LiSo2000,LiSo2001,LoQuMa2005,%
Mas86,%
Mas87a,Mas87b,MaKo94,MaSa92,MaVo88,Qu94,Vor63,Vor67,Vor70,%
WoOl2005,ZiK2005,ZiU81}.

There are many applications to differential equations and
stochastic differential equations, see, e.g., 
\cite{Fl2002,Fl2004,FlMc2000,FlMc2005,Gun98a,Ko96,Ko2000,Ko2001a,%
KoMa97,Mas86,Mas87a,Mas87b,MaKo94,MaSa92,Pap2005,So99a,So99b}.

Applications to the game theory are discussed, e.g., in 
\cite{KoMal97,KoMa97,MaKo94}. There are interesting applications
in biology (bioinformatics), see, e.g., \cite{FiRo93,PaSt2004,Roy92}. Applications
and relations to mathematical morphology are examined the paper 
\cite{DeDo2001} of P.~Del Moral and M.~Doisy and especially in an extended
preprint version of this article. There are many relations and
applications to physics (quantum and classical physics, 
statistical physics, cosmology etc.) see, e.g., Section 6 above and \cite{Ko2000,%
KoMa97,LiMaSh2001,LiMaSh2002,LoQuMa2005,Nu91,Qu97,ChDu87}.

There are important relations and applications to purely mathematical areas.
The so-called tropical combinatorics is discussed in a large survey paper
\cite{Ki2001} of A.~N.~Kirillov, see also \cite{Bu2005,ZiU81}. Interesting
applications of tropical semirings to the traditional representation theory
are presented in \cite{BeFoZe96,BeZe2001,Ki2001}.
Tropical mathematics is closely
related to the very attractive and popular theory of cluster
algebras founded by S.~Fomin and A.~Zelevinsky, see their survey
paper \cite{FoZe2004}. In both cases there are relations with the
traditional theory of representations of Lie groups and related topics.
There are important relations with convex analysis and discrete
convex analysis, see, e.g., \cite{AkGa2003,CoGaQuSi2005,Cu95,DaKo2004,% 
DeSt2003,LiSh2005,MaTi2003,MaSa92,Si97,ZiK77,ZiK79a,ZiK79b}. Some results
on complexity of idempotent/tropical calculations can be found, e.g.,
in \cite{KiRo2005a,LiSo2000,LiSo2001,The2004}.
Interesting applications of tropical algebra to the theory of braids
and the Yang-Baxter mappings (in the sense of \cite{Buh98}) can be found
in \cite{DehDy2002,Dy2002,DyWi2004}.

Many authors examine, explicitly or not, relations and applications of idempotent
mathematics to mathematical economics starting from the 
classical papers of N.~N.~Vorobjev \cite{Vor63,Vor67,Vor70}, see,
e.g., \cite{DaKoMu2001,KoMal97,MaSa92,ZiK76,ZiU81}.

\end{document}

\title[Деквантование Маслова]
{Деквантование Маслова, идемпотентная и тропическая математика:
краткое введение}

\author{Г.~Л.~Литвинов}
\address{Независимый Московский Университет\\
Б.\ Власьевский пер., 11\\
Москва 121002, Россия} 
\email{islc@dol.ru}

\date{}

%\subjclass[2000]{Primary 00A05, 81Q20, 14P99, 51P05, 46S10, 49L99;
% Secondary 06A99, 06F07, 14N10, 81S99.}
\thanks{Эта работа поддержана грантом РФФИ 05--01--00824
и Международным институтом математической физики им. Э.~Шредингера
в Вене (Erwin Schr{\"o}\-din\-ger International Institute for
Mathematical Physics).\\[2ex]
\textit{Записки научных семинаров ПОМИ.  Теория представлений,
  динамические системы~XII}, в печати}

\dedicatory{Анатолию Моисеевичу Вершику с~восхищением и~благодарностью}

\begin{abstract}
Эта статья является кратким введением в идемпотентную и
тропическую математику. Тропическая математика может быть получена
из традиционной математики над числовыми полями как результат
деквантования Маслова, при котором постоянная Планка
 $\hbar$ стремится к нулю, принимая мнимые
значения.

{\it Ключевые слова:} деквантование Маслова, идемпотентная
математика, тропическая математика, идемпотентные полукольца,
идемпотентный анализ, идемпотентный функциональный анализ,
деквантование геометрии.
\end{abstract}

\maketitle

Эта статья практически не содержит строгих формулировок теорем и
их доказательств, она является лишь кратким введением в
деквантование Маслова, идемпотентную и тропическую математику.
Список цитируемой литературы отнюдь не является полным.
Дополнительные ссылки могут быть найдены, например, в электронном
архиве
\[
  \text{\texttt{http://arXiv.org}}
\]
и в работах \cite{BaCoOlQu92,Bu94,%
Ca79,CoGaQu99,CoQu94,Cu79,Cu95,CuMe80,DeDo98,DeDo2001,%
Fl2004,Gla2002,Gol99,Gol2000,Gol2003,GoMi79,GoMi2001,Gun98a,%
Gun98b,KiRo2004,Ki2001,KlMePa2000,Ko2001,KoMa97,LiMa95,LiMa98,%
LiMaSh2001,LiSo2001,Pap2002,ZiU81}. Настоящий обзор является
расширенной версией статьи \cite{Li2005}.

Автор благодарен М.~Акиан (M.\ Akian), Я.\ И.\ Белопольской,
П.~Бу\-т\-ко\-вичу (P.\ Butkovi{\v c}), Г.~Коэну (G.\ Cohen),
С.~Гоберу (S.\ Gaubert), Р.~А.~Ка\-нин\-гхем-Грину (R.\ A.\
Сuninghame-Green), П.~Дел Моралу (P.\ Del Moral), У.~Х.~Флемингу
(W.\ H.\ Fleming), Дж.~С.~Голану (J.\ S.\ Golan), М.~Гондрану (M.\
Gondran), И.~Итенбергу (I.\ Itenberg), Е.~Кацову (Y.\ Katsov), В.\
Н.\ Колокольцову, П.~Лорети (P.\ Loreti), Ю.~И.\ Манину,
Дж.~Маскари (G.\ Mascari), У.~М.~Макэнэни (W.\ M.\ Mac\-Eneaney),
Г.\ Паниной, Е.~Папу (E.\ Pap), М.~Педишини (M.\ Pedicini),
А.~Праду (H.\ Prade), А.~А.~Пухальскому (A.\ A.\ Puhalskii),
Ж.-П.~Квадра (J.-P.\ Quadrat), M.~A.~Ройтбергу, А.~В.~Чуркину,
Г.~Б.~Шпизу, И.~Зингеру (I.\ Singer), и О.~Виро (O.\ Viro) за
полезные добавления и замечания. Автор особенно благодарен
В.~П.~Маслову за всестороннюю помощь и поддержку и
А.~Н.~Соболевскому за большую помощь, включая три рисунка и многое
другое.

\section{Основные идеи}

В основе идемпотентной математики лежит замена обычных
ари\-фметических операций новым набором базовых операций (такими как
максимум или минимум), при этом числовые поля заменяются
идемпотентными полукольцами и полуполями. Типичные примеры~--- так
называемые алгебра макс-плюс $\rmax$ и алгебра мин-плюс $\rmin$.
Пусть $\rset$ - поле вещественных чисел. Тогда ${\rmax}={\rset}
\cup \{-\infty\}$ с операциями $x\oplus y=\max \{x,y\}$ и $x\odot
y=x+y$. Аналогично ${\rmin}={\rset}\cup \{+\infty\}$ с операциями
$\oplus=\min$, $\odot=+$. Новое сложение $\oplus$ является
идемпотентной операцией, то есть $x\oplus x=x$ для всех $x$.

Начиная с классической работы С.~Клини \cite{Kle56}, многие авторы
(C.\ Клини, C.\ Н.\ Н.\ Пандит, Н.~Н.~Воробьев, Б.\ A.\ Карре,
Р.\ A.\ Ка\-нинг\-хем-Грин, K.\ Циммерманн, У.\ Циммерманн, M.\
Гондран, Ф.\ Л.\ Бaккелли, Г.\ Коэн, С.\ Гобер, Г.\ Дж.\ Oлсдер,
Ж.-П.\ Квадра и др.) использовали идемпотентные полукольца и
матрицы над этими полукольцами для решения ряда прикладных задач
дискретной математики и информатики. Современный
\emph{идемпотентный анализ} (или \emph{идемпотентное исчисление},
или \emph{идемпотентная математика}) был разработан В.~П.~Масловым
и его сотрудниками в восьмидесятых годах в Москве, см.,
например, работы
\cite{Mas86,Mas87a,Mas87b,MaKo94,MaSa92,MaVo88,KoMa97}. Некоторые
предварительные результаты сформулировали Э. Хопф и Г. Шоке, см.
работы \cite{Cho55,Ho50}.

Идемпотентную математику можно рассматривать как результат
деквантования традиционной математики над числовыми полями, при
котором постоянная Планка $\hbar$ стремится к нулю, принимая
мнимые значения. Такая точка зрения была представлена
Г.~Л.~Литвиновым и В.~П.~Мас\-ло\-вым
\cite{LiMa95,LiMa96,LiMa98}, см. также
\cite{LiMaSh2001,LiMaSh2002}. Иначе говоря, идемпотентная
математика является асимптотической версией традиционной
математики над полями вещественных и комплексных чисел.

Основную парадигму идемпотентной математики выражает {\it
идемпотентный принцип соответствия}. Этот принцип тесно связан\nopagebreak~со 
знаменитым принципом соответствия Нильса Бора для квантовой
теории. Оказывается, что существует эвристическое соответствие
между рядом важных, интересных и полезных конструкций и
результатов обычной математики над полями и аналогичными
конструкциями и результатами над идемпотентными полуполями и
полукольцами (полуполями и полукольцами с идемпотентными
сложением).

\begin{figure}
  \label{fig:quant}
  \centering
  \includegraphics[width=12cm]{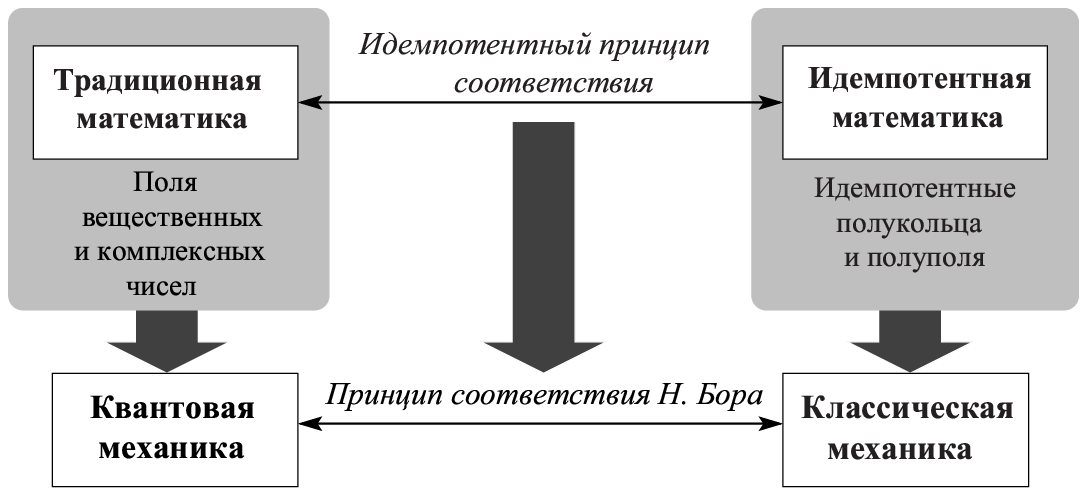}
  \caption{Связь между идемпотентной и традиционной математикой.}
\end{figure}

Систематическое и последовательное использование идемпотентного
принципа соответствия приводит к многообразным результатам, часто
весьма неожиданным. В результате, наряду с традиционной
математикой, возникает ее ``теневая'' идемпотентная версия. Эта
``теневая'' версия так же связана с традиционной математикой, как
классическая физика с физикой квантовой, см. рис.
~\ref{fig:quant}.

Во многих отношениях идемпотентная математика проще традиционной.
Однако переход от традиционных конструкций и результатов к их
идемпотентным аналогам часто является нетривиальным.

%%% {\bf 2. Semirings, semifields, and dequantization.}
\section{Полукольца, полуполя и деквантование}

Пусть на множестве $S$ заданы две алгебраические операции : {\it
сложение} $\oplus$ и {\it умножение} $\odot$. Говорят, что на
множестве $S$ задано {\it полукольцо}, если выполняются следующие
условия:
\begin{itemize}
\item сложение $\oplus$ и умножение $\odot$ ассоциативны;
\item сложение $\oplus$ коммутативно;
\item умножение $\odot$ дистрибутивно относительно сложения $\oplus$: \[x\odot(y\oplus z)=(x\odot y)\oplus(x\odot z)
\quad\text{и}\quad (x\oplus y)\odot z=(x\odot z)\oplus(y\odot z)\]
для любых $x,y,z\in S$.
\end{itemize}
{\it Единицей} полукольца $S$ называется такой элемент
 $\1\in S$,\quad что \quad $\1\odot x=x\odot\1=x$ для всех $x\in
S$. {\it Нулем} полукольца $S$ называется такой элемент $\0\in S$,
что $\0\neq\1$ и $\0\oplus x=x$, $\0\odot x=x\odot \0=\0$ для всех
$x\in S$. Полукольцо $S$ называется {\it идемпотентным
полукольцом}, если $x\oplus x=x$ для всех $x\in S$. Полукольцо $S$
с элементами $\0$ и $\1$ называется {\it полуполем}, если для
любого ненулевого элемента множества $S$ существует обратный
элемент.
Заметим, что диоиды в смысле работ \cite{BaCoOlQu92,GoMi79,%
GoMi2001}, кванталы в  смысле \cite{Ro90,Ro96}, и инклайны в
смысле \cite{KiRo2004} являются примерами идемпотентных полуколец.

\begin{figure}[t]
  \centering
  \includegraphics[height=12cm]{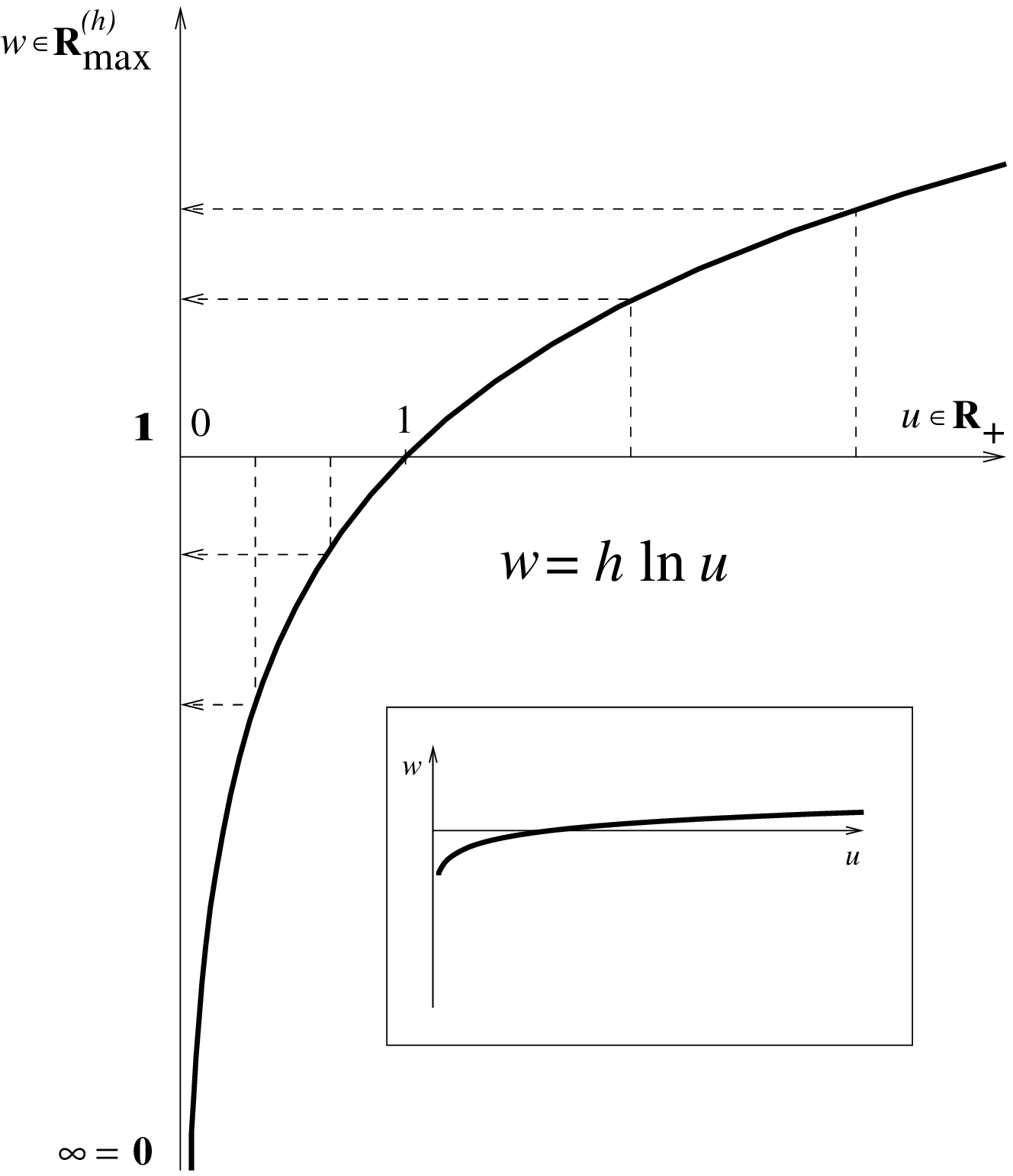}
  \label{fig:dequant}
  \caption{Переход от~\(\rset_+\) к~\(\rset^{(h)}\).  На вставке: то же для малых~\(h\).}
\end{figure}

Рассмотрим поле вещественных чисел $\rset$ и  полуполе всех
неотрицательных вещественных чисел $\rset_+$ (относительно обычных
операций сложения и умножения). Замена переменных $x \mapsto u = h
\ln x$, $h> 0$, задает отображение $\Phi_h \colon \rset_+ \to S =
\rset \cup \{-\infty\}$.  Перенесем операции сложения и умножения
 из $\rset$ в $S$ с помощью отображения  $\Phi_h$, а именно, пусть
$u \oplus_h v = h \ln (\exp (u/h) + \exp(v/h))$, $u \odot v = u +
v$, $\0 = -\infty = \Phi_h(0)$, $\1 = 0 = \Phi_h(1)$, таким
образом $S$ приобретает структуру полукольца $\rset^{(h)}$
изоморфного $\rset_+$; см. рис.~\ref{fig:dequant}.

Несложно проверить, что $u \oplus_h v \to \max \{u,v\}$ при $h \to
0$ и что $S$ образует полукольцо относительно сложения $u \oplus v
= \max \{u,v\}$ и умножения $u \odot v = u + v$ с нулевым
элементом $\0 = -\infty$ и единицей $\1 = 0$. Обозначим это
полукольцо через $\rmax$; оно {\it идемпотентно}, так как $u
\oplus u = u$ для всех элементов. При этом полукольцо
$\rset_{\max}$ является полуполем. Аналогия с процедурой
квантования здесь очевидна, параметр $h$ играет роль постоянной
Планка, поэтому полуполе $\rset_+$  можно рассматривать как
``квантовый'' объект, а полукольцо $\rmax$ может рассматриваться как
результат его ``деквантования''. Аналогичная процедура для $h<0$
дает полукольцо $\rmin = \rset \cup\{+\infty\}$ с операциями
$\oplus = \min$, $\odot = +$; в этом случае $\0 = +\infty$, $\1$ =
0. Полукольца $\rmax$ и $\rset_{\min}$ изоморфны.  Переход к
$\rmax$ или $\rmin$ называется {\it деквантованием Маслова}.
Понятно, что соответствующий переход от $\cset$ или $\rset$ к
$\rset_{\max}$ осуществляется при помощи деквантования Маслова и
отображения $x\mapsto |x|$. Допуская вольность речи, {\it такой
переход в дальнейшем мы также будем называть деквантованием
Маслова}.

Связь с физикой и роль мнимых значений постоянной Планка будут
обсуждаться ниже, в разделе 6, а также в работах
\cite{LiMaSh2001,LiMaSh2002}. Идемпотентное полукольцо $\rset \cup
\{-\infty\} \cup \{+\infty\}$ с операциями $\oplus = \max$, $\odot
= \min$ может быть получено в результате ``вторичного
деквантования''  полей \quad $\cset$, $\rset$ или полуполя
$\rset_+$. Десятки интересных примеров неизоморфных идемпотентных
полуколец, а также стандартные методы получения новых полуколец из
исходных рассматриваются в следующих
работах \cite{CoGaQu2004,Gol99,Gol2003,%
GoMi79,GoMi2001,Gun98a,Gun98b,%
LiMa98,LiMaSh2001}. Так называемое {\it идемпотентное
деквантование} является обобщением деквантования Маслова; это
переход от полей к идемпотентным полуполям и полукольцам в
математических конструкциях и результатах.

Деквантование Маслова связано с хорошо известным логарифмическим
преобразованием, использованным, например, в классических работах
Э.~Шредингера \cite{Sch26} и Э.~Хопфа \cite{Ho50}. Используется
также термин ``преобразование Коула-Хопфа''.
 Идеи Э.~Хопфа получили дальнейшее развитие в
известном методе исчезающей вязкости
и методе вязкостных решений, см., например, работы \cite{BaCa97,%
CaLi97,FlSo93,MaKo94,McC2005,Rou2005,Su95}.

%%%{\bf 3. Tropical semirings and tropical mathematics.}
\section{Терминология: тропические полукольца и тропическая математика}

Термин ``тропические полукольца''  появился в информатике и теории
алгоритмов для обозначения дискретной версии алгебры $\rmax$ или
$\rmin$ и их подалгебр; дискретные полукольца этого типа были
названы ``тропическими'' Домиником  Перрэном в честь бразильского
специалиста по информатике и математике Имре Саймона, в знак
признания его пионерской деятельности в данной области, см.
\cite{Pin98}.

В дальнейшем ситуация и терминология изменились. Для большинства
современных авторов ``тропический'' означает ``над полуполями
$\rset_{\max}$ или $\rset_{\min}$'', а тропические полукольца~--- это
идемпотентные полуполя $\rset_{\max}$ и $\rset_{\min}$. В этом же
смысле часто используются  термины ``макс-плюс'' и ``мин-плюс''. В
настоящее время термин ``тропическая математика'' обычно означает 
``математика над полуполями $\rset_{\max}$ или $\rset_{\min}$'', см.,
например,
\cite{Ale2005,BeFoZe96,BeZe2001,BlYu2005,Dev2004,DeSaSt2004,%
FoGo2003a,FoGo2003b,%
GaMa2005a,GaMa2005b,ItKhSh2003,ItKhSh2003a,ItKhSh2004,ItKhSh2005,Izh2005a,%
Izh2005b,Jo2004,KiRo2005,KiRo2005a,Mi2001,Mi2003,Mi2004,Mi2005,NiSie2004,%
NoYa2002,PaSt2003,PaSt2004,RiStTh2005,Shu2002,Shu2004,So2005,Sp2004,SpSt2003,%
SpSt2004,SpWi2003,Te2004,The2004,Vi2004}. Термины ``тропикализация''
и ``тропификация'' (например, в \cite{Ki2001}) в точности означают
деквантование и квантование в описанном выше смысле. В любом
случае, тропическая математика является естественной и очень
важной частью идемпотентной математики. Многие известные
конструкции и результаты идемпотентной математики  были повторены
в рамках тропической математики (и, особенно, в тропической
линейной алгебре).

Заметим, что в статьях \cite{Vor63,Vor67,Vor70} Н.~Н.~Воробьев
развил некоторую версию идемпотентной линейной алгебры (с важными
приложениями, в том числе для математической экономики), и
предвидел многие аспекты будущей расширенной теории. Для
обозначения идемпотентных полуколец и идемпотентной математики он
использовал термины ``экстремальные алгебры'' и ``экстремальная
математика''. К сожалению, идеи Н.~Н.~Воробьева в свое время не
получили широкой известности, поэтому его терминология не
прижилась и сейчас почти не используется.

%%% {\bf 4. Idempotent algebra and linear algebra.}
\section{Идемпотентная алгебра и линейная алгебра}

Автором первой известной работы по идемпотентной линейной алгебре
является C.~Клини. В его работе \cite{Kle56} рассматриваются системы
линейных алгебраических уравнений над несколько экзотическим
идемпотентным полукольцом всех формальных языков с фиксированным
конечным алфавитом. Однако идеи C.~Клини оказались весьма общими и
универсальными. После этого десятки авторов изучали матрицы с
коэффициентами, принадлежащими идемпотентным полукольцам, а также
соответствующие приложения к дискретной математике, информатике,
языкам программирования, лингвистическим задачам, конечным
автоматам, проблемам оптимизации на графах, теории оптимального
управления, дискретным системам событий и сетям Петри,
стохастическим системам, оценке производительности компьютеров,
вычислительным проблемам и т.д. Эти направления хорошо известны и
широко представлены в литературе, см., например
 \cite{BaCoOlQu92,Bu94,%
Bu2005,Ca79,CoGaQu99,CoQu94,Cu79,Cu95,DuSa92,Gla2002,%
Gol99,Gol2000,Gol2003,GoMi79,GoMi2001,Gun98a,Gun98b,%
KiRo2004,KoMa97,Ko2001,LiMa95,LiMa98,LiMa2005,LiMaE2000,%
LiMaRo2000,LiSo2000,LiSo2001,%
MaKo94,Vor63,Vor67,Vor70,Wa2005,ZiU81}.

Идемпотентная абстрактная  алгебра пока не так хорошо развита,
хотя с другой стороны, с формальной точки зрения теория решеток,
теория упорядоченных групп и полугрупп входят в состав
идемпотентной алгебры. Тем не менее, имеется много интересных
результатов и приложений, см., например,
\cite{Cu79,Cu95,CuMe80,Ka2004,Ro90,Ro96,Shp2000}.

В частности, идемпотентная версия основной теоремы алгебры
сформулирована в \cite{CuMe80,Shp2000} для радикабельных
идемпотентных полуколец (полукольцо $A$ называется {\it
радикабельным}, если уравнение  $x^n = a$ имеет решение $x\in A$
для любого  $a\in A$ и любого положительного целого $n$).
Доказано, что $\rset_{\max}$ и другие радикабельные полуполя
алгебраически замкнуты в естественном смысле \cite{Shp2000}.

В последние годы особое внимание привлекают к себе вопросы
тропической алгебраической геометрии, которые будут рассмотрены
ниже, в разделе 11.

%%%{\bf 5. Idempotent analysis.}
\section{Идемпотентный анализ}

Идемпотентный анализ первоначально был построен В.П. Масловым и
его сотрудниками, а затем развит в большом количестве публикаций
различных авторов. Этому предмету посвящена книга
В.~Н.~Колокольцова и В.~П.~Маслова \cite{KoMa97} (русская версия
\cite{MaKo94} опубликована в 1994 году).

Пусть $S$~--- произвольное полукольцо с идемпотентным сложением
$\oplus$ (которое всегда предполагается коммутативным), умножением
$\odot$, нулем $\0$, и единицей $\1$. Множество $S$ снабжено
 {\it стандартным частичным порядком\/}~$\cle$: по определению,
  $a \cle b$ тогда и только тогда, когда $a \oplus b = b$.
Таким образом, все элементы множества $S$ неотрицательны: $\0
\cle$ $a$ для всех $a \in S$. Благодаря существованию этого
порядка, идемпотентный анализ тесно связан с теорией решеток,
теорией векторных решеток и теорией упорядоченных пространств.
Более того, этот частичный порядок позволяет моделировать основные
топологические понятия и результаты на чисто алгебраическом
уровне; систематическое изложение такого подхода можно найти в
работах \cite{LiMaSh98,LiMaSh99,LiMaSh2001,LiMaSh2002,%
LiSh2002} и \cite{CoGaQu2004}.

Анализ главным образом имеет дело с функциями, чьими значениями
являются числа. Идемпотентным аналогом числовой функции является
отображение $X \to S$, где $X$~--- произвольное множество, а $S$ -
идемпотентное полукольцо. Функции со значениями в $S$ можно
скадывать, умножать друг на друга, а также умножать на элементы
$S$ поточечно.

Идемпотентным аналогом линейного функционального пространства
является множество функций со значениями в  $S$, замкнутое
относительно сложения функций и умножения функций на элементы из
$S$, или $S$-полумодуль. Рассмотрим, например, $S$-полумодуль
$B(X, S)$ всех функций $X \to S$, которые ограничены в смысле
стандартного порядка в $S$.

Пусть $S = \rmax$, тогда идемпотентный аналог интеграла
определяется формулой
$$
I(\varphi) = \int_X^{\oplus} \varphi (x)\, dx = \sup_{x\in X}
\varphi (x),\eqno{(1)}
$$
где $\varphi \in B(X, S)$. Действительно, риманова сумма вида
$\sumlim_i \varphi(x_i) \cdot \sigma_i$ соответствует выражению
$\bigoplus\limits_i \varphi(x_i) \odot \sigma_i = \maxlim_i
\{\varphi(x_i) + \sigma_i\}$, которое сходится к правой части (1)
при $\sigma_i \to 0$. Разумеется, это чисто эвристическое
рассуждение.

Формула~(1) определяет \emph{идемпотентный интеграл} (или
\emph{интеграл Маслова}) не только для функций, принимающих
значения в $\rmax$, но также и в общем случае, при условии, что
любое ограниченное сверху подмножество $S$ имеет точную верхнюю
грань.

\emph{Идемпотентная мера} (или \emph{мера Маслова}) на $X$
определяется как $m_{\psi}(Y) = \suplim_{x \in Y} \psi(x)$, где
$\psi \in B(X,S)$, $Y \subseteq X$. Интеграл по этой мере имеет
вид
$$
   I_{\psi}(\varphi)
    = \int^{\oplus}_X \varphi(x)\, dm_{\psi}
    = \int_X^{\oplus} \varphi(x) \odot \psi(x)\, dx
    = \sup_{x\in X} (\varphi (x) \odot \psi(x)).
    \eqno{(2)}
$$
Очевидно, что если $S = \rmin$, то стандартный порядок $\cle$
является противоположным по отношению к обычному порядку $\le$,
поэтому выражение~(2) примет вид
$$
   \int^{\oplus}_X \varphi(x)\, dm_{\psi}
    = \int_X^{\oplus} \varphi(x) \odot \psi(x)\, dx
    = \inf_{x\in X} (\varphi (x) \odot \psi(x)),
$$
где $\inf$ понимается в смысле обычного порядка $\le$.

Отметим, что так называемый псевдо-анализ (см., например, обзор
\cite{Pap2005}) связан со специальной частью идемпотентного
анализа; тем не менее, вообще говоря, этот псевдо-анализ  выходит
за рамки идемпотентной математики. Некоторые обобщения мер Маслова
обсуждаются в работах \cite{KlMePa2000,Pap2002}.

%% {\bf 6. The superposition principle and linear problems.}
\section{Принцип суперпозиции и линейные задачи}

Основные уравнения квантовой теории линейны. Это утверждение
составляет принцип суперпозиции в квантовой механике. Уравнение
Гамильтона-Якоби, которое является основным уравнением
классической механики, не является линейным в обычном смысле. Тем
не менее, оно линейно над полукольцами $\rmax$ и $\rmin$. Похожим
образом различные версии уравнения Беллмана, основного уравнения
теории оптимизации, линейны над подходящими полукольцами. В этом
состоит
идемпотентный принцип суперпозиции В.П.Маслова (см. \cite{Mas86,Mas87a,Mas87b,MaKo94,%
MaSa92}). Например, конечномерное стационарное уравнение Беллмана
может быть записано в виде  $X = H \odot X \oplus F$, где $X$,
$H$, $F$~--- матрицы с коэффициентами в идемпотентном полукольце
$S$, а неизвестная матрица  $X$ определяется через $H$ и $F$
\cite{Ca71}. В частности, стандартные проблемы динамического
программирования и задача поиска кратчайшего пути отвечают случаям
$S = \rmax$ и $S =\rmin$ соответственно. В работе \cite{Ca71}
оптимизационные алгоритмы на конечном графе поставлены в
соответствие стандартным методам решения систем линейных уравнений
этого типа (над полукольцами). В частности, алгоритму Беллмана
поиска кратчайшего пути соответствует алгоритм Якоби, алгоритму
Форда отвечает итерационный метод Гаусса-Зейделя и т.д.

Линейность над $\rmin$ и $\rmax$ уравнения Гамильтона-Якоби
(которое является результатом деквантования Маслова уравнения
Шредингера) тесно связана с обычной линейностью уравнения
Шредингера и может быть выведена из этой линейности. Таким
образом, можно заимствовать стандартные идеи и методы из линейного
анализа и применять их в новой области.

Рассмотрим классическую динамическую систему, описываемую
гамильтонианом
$$
   H = H(p,x) = \sum_{i=1}^N \frac{p^2_i}{2m_i} + V(x),
$$
где $x = (x_1, \dots, x_N)$~--- обобщенные координаты, $p = (p_1,
\dots, p_N)$~--- обобщенные импульсы, $m_i$~--- массы, а $V(x)$
-потенциал. В этом случае лагранжиан системы $L(x, \dot x, t)$
запишется в виде
$$
   L(x, \dot x, t)
        = \sum^N_{i=1} m_i \frac{\dot x_i^2}2~--- V(x),
$$
где $\dot x = (\dot x_1, \dots, \dot x_N)$, $\dot x_i = dx_i /
dt$. При этом функция значения функционала действия $S(x,t)$ имеет
вид
$$
   S = \int^t_{t_0} L(x(t), \dot x(t), t)\, dt,
$$
где интегрирование ведется вдоль реальной траектории системы.
Классические уравнения движения выводятся из условия минимума
функционала действия (принцип Гамильтона или принцип наименьшего
действия).

Для фиксированных моментов времени $t$ и $t_0$ и произвольной
траектории $x(t)$, функционал действия $S=S(x(t))$ может быть
рассмотрен как функция, переводящая множество кривых (траекторий)
во множество вещественных чисел, которое, в свою очередь, может
быть интерпретировано как $\rmin$. В этом случае минимум
функционала действия выражается интегралом Маслова по множеству
траекторий, т.е. как идемпотентный аналог эвклидовой версии
фейнмановского континуального интеграла.

Минимум функционала действия соответствует максимуму функции
 $e^{-S}$, т.е. идемпотентному интегралу $\int^{\oplus}_{\{paths\}} e^{-S(x(t))} D\{x(t)\}$
по отношению к алгебре $\rset_{\max}$. Таким образом, классический
принцип наименьшего действия может рассматриваться как
идемпотентная версия фейнмановской формулировки квантовой
механики. Представление решения уравнения Шредингера в виде
континуального интеграла Фейнмана соответствует формуле
Лакса-Олейник для уравнения Гамиль\-тона-Якоби.

Поскольку $\partial S/\partial x_i = p_i$, $\partial S/\partial t
= -H(p,x)$, справедливо следующие уравнение Гамильтона-Якоби:
$$
   \pd{S}{t} + H \left(\pd{S}{x_i}, x_i\right)= 0.\eqno{(3)}
$$
Квантование (см., например,\cite{FeHi65}) приводит к уравнению
Шредингера
$$
   -\frac{\hbar}i \pd{\psi}{t}= \widehat H \psi = H(\hat p_i, \hat x_i)\psi,
   \eqno{(4)}
$$
где $\psi = \psi(x,t)$~--- волновая функция, т.е. зависящий от
времени элемент гильбертового пространства $L^2(\rset^N)$, а
$\widehat H$~--- оператор энергии, полученный из классической
функции Гамильтона заменой переменных $p_i$ операторами импульсов
 $\widehat p_i = {\hbar\over i}{\partial \over
\partial x_i}$,
и переменных $x_i$~--- операторами координат:
 $\widehat x_i \colon \psi \mapsto
x_i\psi$. Уравнение Шредингера линейно в обычном смысле (квантовый
принцип суперпозиции).

В стандартной процедуре предельного перехода от уравнения Шредингера
к уравнению Гамильтона-Якоби используется следующее представление
для волновой функции в виде $\psi(x,t) = a(x,t)\cdot\allowbreak e^{iS(x,t)/\hbar}$
и выделение главного члена при $\hbar \to 0$ (квазиклассический
предел).

Вместо этого мы заменим постоянную Планка $\hbar$ мнимой величиной
$h = i\hbar$, выбирая $h> 0$. Тогда уравнение Шредингера (4)
перейдет в обобщенное уравнение теплопроводности:
$$
   h\pd{u}{t} = H\left(-h\frac{\partial}{\partial x_i}, \hat x_i\right) u,
   \eqno{(5)}
$$
в котором вещественнозначная функция $u$ соответствует волновой
функции~$\psi$ (точнее, ее модулю). Похожая идея (переход к мнимому времени)
используется в эвклидовой квантовой теории поля (см., например,
\cite{Ne73}). Напомним, что время и энергия~--- дуальные величины.

Линейность уравнения~(4) ведет к линейности уравнения~(5).
Поэтому, если $u_1$ и $u_2$~--- решения уравнения~(5), то их
линейная комбинация
$$
   u = \lambda_1 u_1 + \lambda_2 u_2\eqno{(6)}
$$
также является решением уравнения (5). Положим $S = h \ln u$ или
$u = e^{S/h}$, как это было сделано выше в разделе 2. Несложно
проверить, что после такой замены уравнение~(5) примет вид
$$
   \pd{S}{t}= V(x) + \sum^N_{i=1} \frac1{2m_i}\left(\pd{S}{x_i}\right)^2
        + h\sum^n_{i=1}\frac1{2m_i}\frac{\partial^2 S}{\partial x^2_i}.
   \eqno{(7)}
$$
Таким образом, мы перешли от (4) к (7) при помощи замены
переменных $\psi = e^{S/h}$. При этом $|\psi| = e^{\Re S/h}$ , где
$\Re S$~--- вещественная часть $S$. Будем рассматривать теперь $S$
как вещественную переменную. Уравнение (7) нелинейно в обычном
смысле. Однако, если $S_1$ и $S_2$~--- его решения, то функция
$$
   S = \lambda_1 \odot S_1 \opl_h \lambda_2\odot S_2,
$$
полученная из (6) заменой $S = h \ln u$, также является решением
уравнения (7). В этом случае обобщенное умножение $\odot$
совпадает с обычным сложением, а обобщенное сложение $\opl_h$
определяется выбранной заменой переменных. При $h \to 0$, мы
получим операции идемпотентного полукольца  $\rmax$, то есть,
$\oplus = \max$ и $\odot = +$, а уравнение (7) перейдет в
уравнение Гамильтона-Якоби (3), так как третье слагаемое в правой
части уравнения~(7) исчезает.

Естественно рассматривать  предельную функцию $S = \lambda_1 \odot
S_1 \oplus \lambda_2 \odot S_2$ как решение уравнения
Гамильтона-Якоби и ожидать, что это уравнение можно рассматривать
как линейное над $\rmax$. Подобная аргументация (на эвристическом
уровне) может быть применена и к уравнениям более общего вида.
Строгое обоснование линейности этих уравнений над полукольцами
приводится в \cite{FlMc2000,KoMa97,MaKo94,MaSa92} а также в
\cite{Mas87a}. Отметим, что если $h$ заменить на $-h$, то в
результате мы получим линейность уравнения Гамильтона-Якоби над
полукольцом $\rmin$.

Идемпотентный принцип суперпозиции показывает, что ряд важных
нелинейных (в обычном смысле) задач может быть сведен к задачам,
линейным над идемпотентными полукольцами. Таким образом, линейный
идемпотентный функциональный анализ (см. ниже) является
естественным инструментом для исследования нелинейных
бесконечномерных задач, обладающих этим свойством.

\section{Свертка и преобразование Фурье--Лежандра}

Пусть $G$~--- некоторая группа. Тогда пространство $B(X,
\rset_{\max})$ всех ограниченных функций $G\to\rset_{\max}$ (см.
выше) является идемпотентным полукольцом  относительно следущего
аналога обычной операции свертки $\circledast$:
$$
   (\varphi(x)\circledast\psi)(g)=
        \int_G^{\oplus} \varphi (x)\odot\psi(x^{-1}\cdot g)\, dx=
\sup_{x\in G}(\varphi(x)+\psi(x^{-1}\cdot g)).
$$
Разумеется, можно также использовать и другие функциональные
пространства (и другие базовые полукольца вместо $\rset_{\max}$).
В работах \cite{KoMa97,MaKo94} групповые полукольца этого типа
называются {\it сверточными полукольцами}.

Пусть $G=\rset^n$, где $\rset^n$ понимается как топологическая
группа относительно векторного сложения. Обычное преобразование
Фурье-Лапласа определяется формулой
$$
   \varphi(x) \mapsto \tilde{\varphi}(\xi)
        = \int_G e^{i\xi \cdot x} \varphi (x)\, dx,\eqno{(8)}
$$
где $e^{i\xi \cdot x}$ является характером группы $G$ и решением
следующего функционального уравнения
$$
   f(x + y) = f(x)f(y).
$$
Идемпотентный аналог этого уравнения имеет вид
$$
   f(x + y) = f(x) \odot f(y) = f(x) + f(y),
$$
поэтому ``непрерывные идемпотентные характеры'' являются линейными
функциями вида $x \mapsto \xi \cdot x = \xi_1 x_1 + \dots + \xi_n
x_n$. В результате преобразование~(8) переходит в преобразование
$$
   \varphi(x) \mapsto \tilde{\varphi}(\xi)
        = \int_G^\oplus \xi \cdot x \odot \varphi (x)\, dx
   = \sup_{x\in G} (\xi \cdot x + \varphi (x)).\eqno{(9)}
$$
Преобразование~(9) есть не что иное, как {\it преобразование
Лежандра\/} (с точностью до обозначений) \cite{Mas87b}.
Преобразования этого вида устанавливают связь между гамильтоновым
и лагранжевым формализмом в классической механике. Преобразование
Лежандра порождает идемпотентную версию гармонического анализа в
пространстве выпуклых функций (подробнее см., например,
\cite{MaTi2003}).

Разумеется, описанную конструкцию можно обобщить на различные
классы групп и полуколец. Преобразования указанного типа переводят
обобщенную операцию свертки в поточечное произведение
(обобщенное), сохраняя при этом ряд важных свойств обычного
преобразования Фурье. Для случая полукольца  множеств Парето,
соответствующая версия преобразования Фурье сводит задачу
многокритериальной оптимизации к семейству однокритериальных задач
\cite{SaTa92}.

Примеры, обсуждаемые в данном разделе, могут рассматриваться как
фрагменты идемпотентной  теории представлений, см., например,
\cite{LiMaSh2002}. В частности, ``идемпотентные''  представления
групп можно трактовать как представления соответствующих
сверточных полуколец (т.е. идемпотентных групповых полуколец) в
идемпотентных полумодулях.

%% {\bf 8. Correspondence to stochastics and a duality between
%% probability and optimization.}
\section{Приложения к стохастике и двойственность между теорией вероятностей и теорией оптимизации}

Идемпотентные меры Маслова являются неотрицательными (в смысле
стандартного порядка), как и вероятностные меры. Аналогия между
идемпотентными и вероятностными мерами определяет важную
взаимосвязь между теорией оптимизации и теорией вероятностей. К
настоящему времени исследован целый ряд идемпотентных аналогов
вероятностных конструкций, например макс-плюс мартингалы,
макс-плюс стохастические дифференциальные уравнения и др.
Например, полученные результаты позволяют перенести мощные
стохастические методы в теорию оптимизации. Это было отмечено и
исследовано многими авторами (Г.~Салю, П.~Дел Морал, M.~Акиан,
Ж.-П.~Квадра, В.~П.~Maслов, В.~Н.~Колокольцов, П.~Бернар,
У.~A.~Флеминг, У.~M.~Макэнэни, A.~A.~Пухальский и др.), см
обзорную статью У.~А.~Флеминга и У.~М.~Макэнэни \cite{FlMc2005}
и~\cite{Ak99,AkQuVi98,Be94,CoQu94,De97,De2004,DeDo98,%
DeDo2001,%
Fl2002,Fl2004,FlMc2000,Gun98a,MaKo94,Puh2001,Qu90,Qu94}.
По поводу связей и приложений к теории больших уклонений см. \cite{Ak99,%
De97,De2004,DeDo98,DeDo2001,Puh2003}, и особенно книгу
A.~A.~Пухальского \cite{Puh2001}.

%% {\bf 9. Idempotent functional analysis}.
\section{Идемпотентный функциональный анализ}

Для многих конструкций и результатов традиционной математики можно указать
интересные идемпотентные аналоги. В частности, это относится к
основным конструкциям и теоремам функционального анализа.
Идемпотентный функциональный анализ является абстрактной версией
идемпотентного анализа.  Для простоты изложения положим
 $S=\rmax$, и пусть $X$~--- произвольное множество.
 Идемпотентное интегрирование  задается формулой (1), написанной выше.
Функционал $I(\varphi)$ линеен над $S$, и его значения
соответствуют предельным значениям соответствующих аналогов сумм
Лебега (или Римана). Идемпотентное скалярное произведение двух
функций $\varphi$ и $\psi$ определяется формулой
$$
\langle\varphi,\psi\rangle = \int^{\oplus}_X \varphi(x)\odot\psi(x)\, dx =
\sup_{x\in X}(\varphi(x)\odot\psi(x)).
$$
Поэтому естественно построить идемпотентные аналоги интегральных
операторов в виде
$$
\varphi(y) \mapsto (K\varphi)(x) = \int^{\oplus}_Y K(x,y)\odot
\varphi(y)\, dy = \sup_{y\in Y}\{K(x,y)+\varphi(y)\},\eqno(10)
$$
где $\varphi(y)$~--- элемент пространства функций, определенных на
множестве $Y$, а $K(x,y)$~--- функция, заданная на $X\times Y$, и
принимающая значения в $S$. Как известно, выражения подобного типа
стандартны для теории оптимизации.\medskip

Напомним, что описанные выше определения и конструкции могут быть
успешно обобщены на случай идемпотентных полуколец, которые
условно полны в смысле стандартного порядка. Используя интеграл
Маслова, можно сконструировать различные функциональные
пространства, а также идемпотентную версию теории обобщенных
функций (распределений). Для некоторых конкретных идемпотентных
функциональных пространств доказано, что каждый ``хороший'' линейный
(в идемпотентном смысле) оператор может быть представлен в виде
(10). Это утверждение является идемпотентной версией теоремы
Л.~Шварца о ядре; результаты такого рода были получены
В.~Н.~Колокольцовым, П.~С.~Дудниковым, и С.~Н.~Самборским,
И.~Зин\-гером, М.~А.~Шубиным и другими авторами, см.,например,
\cite{DuSa92,KoMa97,MaKo94,MaSa92,Si2003}. Таким образом, каждый
``хороший'' линейный функционал представим в виде
$\varphi\mapsto\langle\varphi,\psi\rangle$, где $\langle,\rangle$
- идемпотентное скалярное произведение.\medskip

В рамках идемпотентного функционального анализа  результаты этого
типа могут быть доказаны в очень общей ситуации.
В работах \cite{LiMaSh98,%
LiMaSh99,LiMaSh2001,LiMaSh2002,LiSh2002} развита алгебраическая
версия идемпотентного функционального анализа. Это означает, что
основные топологические понятия и результаты моделируются в чисто
алгебраических терминах. Построенная теория охватывает предмет,
начиная с основных понятий и результатов (например, идемпотентных
аналогов знаменитых теорем Хана-Банаха, Банаха-Штейнгауза, Рисса и Рисса-Фишера) и
включая идемпотентные аналоги определений и результатов
А.Гротендика (А.Grothendieck) по топологическим тензорным
произведениям, ядерным пространствам и операторам. Сформулирована
абстрактная версия теоремы о ядре. Отметим, что переход от обычной
теории к идемпотентному функциональному анализу может быть весьма
нетривиальным. Так, например, существует много неизоморфных
идемпотентных гильбертовых пространств. Важные результаты
идемпотентного функционального анализа (теоремы о двойственности и
отделимости) недавно опубликовали Г.~Коэн, С.~Гобер, и
Ж.-П.~Квадра \cite{CoGaQu2004}; конечномерная версия теоремы об
отделимости опубликована в \cite{ZiK77}. Последнее время
идемпотентный функциональный анализ привлекает повышенное
внимание, см., например, \cite{AkGa2003,AkGaKo2005,AkGaWa2004,%
AkGaWa2005,CoGaQuSi2005,GoMi2001,LiMa2005,LiSh2005,%
Rou2005,Si2003,So99a,Wal2005} и работы, цитированные выше.

%%%\section{10. The dequantization transform and generalization of the
%%%Newton polytopes}
\section{Преобразование деквантования и обобщение многогранников Ньютона}

В этом разделе мы кратко обсудим результаты, опубликованные в
работе \cite{LiSh2005}. Для функций, определенных на $\cset^n$
почти всюду, можно построить преобразование деквантования
$f\to\hat f$, порожденное деквантованием Маслова. Для
полиномиальной функции $f$ субдифференциал $\partial{\hat f}$
функции $\hat f$ совпадает с многогранником Ньютона полинома $f$.
Для полукольца полиномов с неотрицательными коэффициентами
преобразование деквантования является гомоморфизмом этого
полукольца в идемпотентное полукольцо выпуклых многогранников
относительно хорошо известных операций Минковского.
Эти результаты
могут быть обобщены на широкий класс функций и выпуклых множеств.
\smallskip

\noindent\textbf{10.1. Преобразование деквантования.}
Рассмотрим топологическое пространство $X$. Для функций $f(x)$,
определенных на $X$, будем говорить, что некоторое свойство
выполняется {\it почти всюду} (п.в.), если оно  справедливо для
всех элементов $x$ из открытого всюду плотного подмножества в $X$.
Выберем в качестве $X$ множество $\cset^n$ или $\rset^n$;
обозначим через $\rset^n_+$ множество $x=\{\,(x_1, \dots, x_n)\in
X \mid x_i\geq 0$, где $i = 1, 2, \dots, n\}$. Для $x= (x_1,
\dots, x_n)\in X$ положим $\exp(x) = (\exp(x_1), \dots,
\exp(x_n))$; таким образом, если $x\in\rset^n$, то $\exp(x)\in
\rset^n_+$.

Обозначим через $\maF(\cset^n)$ множество всех функций,
непрерывных на открытом всюду плотном подмножестве $U\subset
\cset^n$, причем $U\supset \rset^n_+$. Во всех ниже приведенных
примерах мы рассмотрим даже более регулярные функции, которые
голоморфны в~$U$.
 Очевидно, что $\maF(\cset^n)$~--- кольцо
(и алгебра над $\cset$) по отношению к обычному сложению и
умножению функций.

Для $f\in \maF(\cset^n)$ определим функцию $\hat f_h$ по следующей
формуле:
$$
\hat f_h(x) = h \log|f(\exp(x/h))|,\eqno(11)
$$
где $h$ (малый) вещественный параметр, а $x\in\rset^n$. Положим
$$
\hat f(x) = \lim_{h\to 0} \hat f_h (x),\eqno(12)
$$
если предел в правой части равенства (12) существует почти всюду.
Назовем функцию $\hat f(x)$ {\it деквантованием} функции $f(x)$, а
преобразование $f(x)\mapsto \hat f(x)$~--- {\it преобразованием
деквантования}. По построению $\hat f_h(x)$ и $\hat f(x)$ могут
быть рассмотрены как функции, принимающие значения в
$\rset_{\max}$. Отметим, что на самом деле $\hat f_h(x)$ и $\hat
f(x)$ зависят только от сужения функции $f$ на подмножество
$\rset_+^n $; поэтому фактически преобразование деквантования
задается для функций, определенных только на $\rset^n_+$. Понятно,
что преобразование деквантования порождено деквантованием Маслова
и отображением $x\mapsto |x|$. Разумеется, аналогичные определения
могут быть даны и для функций, заданных на $\rset^n$ и
$\rset_+^n$.

Обозначим через $V$ множество $\rset^n$, рассматриваемое как
линейное эвклидово пространство  (со скалярным произведением $(x,
y) = x_1y_1+ x_2y_2 +\dots + x_ny_n$) и положим $V_+ = \rset_+^n$.
Назовем функцию  $f\in \maF(\cset^n)$  {\it деквантуемой}, если ее
деквантование $\hat f(x)$ существует (и определено на открытом
всюду плотном подмножестве в $V$). Обозначим через $\maD
(\cset^n)$ множество всех деквантуемых функций и пусть
$\widehat{\maD}(V)$ означает множество $\{\,\hat f \mid f\in
\maD(\cset^n)\,\}$. Напомним, что функции из $\maD(\cset^n)$ (и
$\widehat{\maD}(V)$) определены почти всюду и равенство $f=g$
означает, что $f(x) = g(x)$ для любого $x$ из открытого всюду
плотного подмножества в $\cset^n$ (соответственно в $V$).
Обозначим через $\maD_+(\cset^n)$ множество всех функций $f\in
\maD(\cset^n)$ таких, что $f(x_1, \dots, x_n)\geq 0$, если
$x_i\geq 0$ для $i= 1,\dots, n$; тогда
 $f\in \maD_+(\cset^n)$, если сужение $f$ на $V_+ = \rset_+^n$ является неотрицательной функцией.
 Обозначим полученный в результате деквантования
 образ множества  $\maD_+(\cset^n)$ через $\widehat{\maD}_+(V)$.
Говорят, что функции $f, g\in \maD(\cset^n)$ находятся в {\it
общем положении}, если $\hat f (x) \neq \widehat g(x)$ для
элементов $x$, пробегающих открытое всюду плотное множество в $V$.

Для функций $f, g \in \maD(\cset^n)$ и ненулевой константы $c$,
верны следующие равенства:

\begin{enumerate}
\item[1)] $\widehat{fg} = \hat f + \widehat g$;
\item[2)] $|\hat f| = \hat f$; $\widehat{cf} = f$; $\widehat c =0$;
\item[3)] $(\widehat{f+g})(x) = \max\{\hat f(x), \widehat g(x)\}$ п.в., если $f$ и $g$
неотрицательны на $V_+$ (то есть $f, g \in \maD_+(\cset^n)$), или
$f$ и $g$~--- функции общего положения.
\end{enumerate}
Левые части этих уравнений определены автоматически.

Множество $\maD_+(\cset^n)$ имеет естественную структуру
полукольца относительно обычного сложения  и умножения функций,
принимающих значение в $\cset$. Множество $\widehat{\maD}_+(V)$
имеет естественную структуру идемпотентного полукольца
относительно операций $(f\oplus g)(x) = \max \{ f(x), g(x)\}$,
$(f\odot g)(x) = f(x) + g(x)$; элементы из
 $\widehat{\maD}_+(V)$ можно рассматривать как функции,
 принимающие значения в $\rset_{\max}$.  Преобразование деквантования порождает гомоморфизм из
 $\maD_+(\cset^n)$ в $\widehat{\maD}_+(V)$.
\smallskip

\noindent\textbf{10.2. Простые функции.}
Для любого ненулевого числа $a\in\cset$ и любого вектора
 $d = (d_1, \dots, d_n)\in V = \rset^n$
положим $m_{a,d}(x) = a \prod_{i=1}^n x_i^{d_i}$; функции такого
вида будем называть {\it обобщенными мономами}. Обобщенные мономы
определены п.в. на  $\cset^n$ и на  $V_+$, но не на $V$, исключая
случай, когда $d_i$ являются целыми или подходящими рациональными
числами. Назовем функцию $f$ {\it обобщенным полиномом}, если она
представляет собой конечную сумму линейно независимых обобщенных
мономов.
  Например, полиномы Лорана являются обобщенными полиномами.

Как обычно, для $x, y\in V$ положим $(x,y) = x_1y_1 + \dots +
x_ny_n$. Легко показать, что если $f$~--- обобщенный моном
$m_{a,d}(x)$, то $\hat f$~--- линейная функция $x\mapsto (d,x)$.
Если $f$~--- обобщенный полином, то $\hat f$~--- сублинейная функция.

Напомним, что вещественная функция $p$, определенная на $V =
\rset^n$, {\it сублинейна}, если $p = \sup_{\alpha} p_{\alpha}$,
где $\{p_{\alpha}\}$~--- семейство линейных функций. Сублинейные
функции, определенные всюду на $V=\rset^n$, выпуклы;
следовательно, они непрерывны. В дальнейшем мы будем рассматривать
только такие сублинейные функции. Пусть $p$~--- непрерывная функция
на $V$, тогда $p$ сублинейна в точности тогда, когда

1) $p(x+ y) \leq p(x) + p(y)$ для всех $x, y \in V$;

2) $p(cx) = cp(x)$ для всех $x\in V$, $c\in \rset_+$.

Если $p_1$, $p_2$~--- сублинейные функции, то $p_1 +p_2$ также
является сублинейной функцией.

Назовем функцию $f \in \maF(\cset^n)$ {\it простой}, если ее
деквантование $\hat f$ существует и п.в.\ совпадает с некоторой
сублинейной функцией. Допуская вольность речи, будем обозначать эту
(однозначно определенную всюду на $V$) сублинейную функцию тем же
символом~$\hat f$.

Напомним, что простые функции $f$ и $g$ являются функциями 
{\it в~общем положении}, если $\hat f(x) \neq \widehat g(x)$ для всех $x$
из открытого всюду плотного множества в $V$. В частности,
обобщенные мономы находятся в общем положении тогда и только
тогда, когда они линейно независимы.

Обозначим через $\mathit{Sim}(\cset^n)$~--- множество всех простых
функций, определенных на $V$, и обозначим через \quad
$\mathit{Sim}_+(\cset^n)$\quad множество $\mathit{Sim}(\cset^n)
\cap \maD_+(\cset^n)$. Через $\mathit{Sbl}(V)$ обозначим множество
всех
  (непрерывных) сублинейных функций, определенных на $V =
  \rset^n$,
а через
 $\mathit{Sbl}_+(V)$ обозначим образ $\widehat{\mathit{Sim}_+}(\cset^n)$
 множества
$\mathit{Sim}_+(\cset^n)$ при преобразовании деквантования.

Множество \quad $\mathit{Sim}_+(\cset^n)$ \quad является
подполукольцом в $\maD_+(\cset^n)$, и $\mathit{Sbl}_+(V)$ является
идемпотентным подполукольцом в $\widehat{\maD_+}(V)$.
Преобразование деквантования порождает \quad эпиморфизм \quad
полукольца $\mathit{Sim}_+(\cset^n)$ на $\mathit{Sbl}_+(V)$.
Множество $\mathit{Sbl}(V)$ является идемпотентным полукольцом
относительно операций $(f\oplus g)(x) = \max \{ f(x), g(x)\}$,
$(f\odot g)(x) = f(x) + g(x)$.

Конечно, полиномы и обобщенные полиномы являются простыми
функциями.

Будем говорить, что функции $f, g\in\maD(V)$ {\it асимптотически
эквивалентны}, если $\hat f = \widehat g$; простую функцию $f$
назовем {\it асимптотическим мономом}, если $\hat f$~--- линейная
функция. Простая функция $f$ будет называться {\it асимптотическим
полиномом}, если $\hat f$ является суммой конечного числа
неэквивалентных асимптотических мономов. Каждый асимптотический
полином является простой функцией.

{\bf Пример.} Примерами асимптотических полиномов служат
обобщенные полиномы, логарифмические функции от (обобщенных)
полиномов, произведения полиномов и логарифмических функций. Это
следует из наших определений и формулы~(11).
\smallskip

\noindent\textbf{10.3. Субдифференциалы сублинейных функций и множества Ньютона для простых функций.}
Хорошо известно, что совокупность всех выпуклых компактных
подмножеств в $\rset^n$ образует идемпотентное полукольцо
$\mathcal{S}$ относительно операций Минковского: для $A, B \in
\mathcal{S}$ сумма $A\oplus B$ является выпуклой оболочкой
объединения $A\cup B$; произведение $A\odot B$ определяется
следующим образом: $A\odot B = \{\, x\mid x = a+b, \text{ где
$a\in A, b\in B$}\,\}$. На самом деле множество $\mathcal{S}$
является идемпотентным линейным пространством над $\rset_{\max}$
(см., например, работу \cite{LiMaSh2001}). Разумеется,
многогранники Ньютона в пространстве $V$ образуют подполукольцо
$\mathcal{N}$ в $\mathcal{S}$.

Воспользуемся некоторыми элементарными результатами из выпуклого
анализа. Эти результаты можно найти, например, в книге
\cite{MaTi2003}.

Для любой функции $p\in \mathit{Sbl}(V)$ положим
$$
\partial p = \{\, v\in V\mid (v, x) \le p(x)\ \forall x\in V\,\}.
$$

Из выпуклого анализа хорошо известно, что для любой сублинейной
функции $p$ множество $\partial p$ есть в точности {\it
субдифференциал} от $p$ в нулевой точке.
 Следующее утверждение также хорошо известно из выпуклого анализа:

\begin{prop}
Пусть $p_1,p_2\in \mathit{Sbl}(V)$, тогда
\begin{center}
$\partial (p_1+p_2) = \partial p_1\odot\partial p_2 =$

$ \{\, v\in V\mid v = v_1+v_2, \text{ где $v_1\in \partial p_1,
v_2\in \partial p_2$}\,\}$;

$\partial (\max\{p_1(x), p_2(x)\}) = \partial p_1\oplus\partial p_2$.
\end{center}

Пусть $p\in \mathit{Sbl}(V)$.  Тогда $\partial p$~--- непустое
выпуклое компактное подмножество в $V$.
\end{prop}

\begin{cor}
Отображение $p\mapsto \partial p$ является гомоморфизмом
идемпотентного полукольца
 $\mathit{Sbl}(V)$ в идемпотентное полукольцо $\mathcal{S}$ всех выпуклых компактных подмножеств множества $V$.
\end{cor}

Для любой простой функции $f\in \mathit{Sim}(\cset^n)$ обозначим
через $N(f)$ множество $\partial(\hat f)$. Назовем $N(f)$ {\it
множеством Ньютона} для функции $f$. Из приведенного утверждения
следует, что для любой простой функции $f$ ее множество Ньютона
$N(f)$ является непустым выпуклым компактным подмножеством в $V$.

\begin{thm}
Пусть $f$ и $g$~--- простые функции. Тогда
\begin{enumerate}
\item[1)] $N(fg) = N(f)\odot N(g) = \{\, v\in V\mid v = v_1 +v_2$
, где  $v_1 \in N(f), v_2 \in N(g)\}$;
\item[2)] $N(f+g) = N(f)\oplus N(g)$, если $f_1$ и $f_2$~--- функции в общем положении или
$f_1, f_2 \in \mathit{Sim}_+(\cset^n)$ {\rm (}напомним, что
$N(f)\oplus N(g)$ есть выпуклая оболочка множества $N(f)\cup
N(g)${\rm )}.
\end{enumerate}
\end{thm}

\begin{cor}
Преобразование $f\mapsto N(f)$ порождает гомоморфизм из
$\mathit{Sim}_+(\cset^n)$ в $\mathcal{S}$.
\end{cor}

\begin{prop}
Пусть $f = m_{a,d}(x) = a \prod^n_{i=1} x_i^{d_i}$~--- (обобщенный)
моном; здесь $d= (d_1, \dots, d_n) \in V= \rset^n$ и $a$ -
ненулевое комплексное число. Тогда $N(f) = \{ d\}$.
\end{prop}

\begin{cor}
Пусть $f = \sum_{d\in D} m_{a_d,d}$~--- (обобщенный) полином. Тогда
$N(f)$~--- многогранник $\oplus_{d\in D}\{d\}$, т.е. линейная
оболочка конечного множества $D$.
\end{cor}

Этот вывод следует из Теоремы и  Предложения 2. В этом случае
 $N(f)$~--- хорошо известный классический многогранник Ньютона полинома $f$.
Поэтому следующее следствие очевидно:

\begin{cor}
Пусть $f$~--- обобщенный или асимптотический полином. Тогда его
множество Ньютона
 $N(f)$ является выпуклым многогранником.
\end{cor}

{\bf Пример.} Рассмотрим одномерный случай, т.е. $V = \rset$ и
положим
 $f_1 = a_nx^n + a_{n-1}x^{n-1} + \dots + a_0$ и $f_2 = b_mx^m + b_{m-1}
 x^{m-1} + \dots + b_0$, где $a_n\neq 0$, $b_m\neq 0$, $a_0 \neq 0$,
 $b_0 \neq 0$. Тогда $N(f_1)$ совпадает с  отрезком $[0, n]$ и $N(f_2)$ -
 с отрезком $[0, m]$. Поэтому преобразование $f\mapsto N(f)$
 соответствует преобразованию $f\mapsto \deg (f)$, где $\deg(f)$~--- степень полинома $f$. При этом в соответствии с
 Теоремой $\deg(fg) = \deg f + \deg g$ и
 $\deg (f+g) = \max \{\deg f, \deg g\} = \max \{n, m\}$, если $a_i\geq 0$,
 $b_i\geq 0$ или если $f$ и $g$~--- находятся в общем положении.

%%{\bf 11. Dequantization of geometry.}
\section{Деквантование геометрии}

Идемпотентная версия  вещественной алгебраической геометрии была
открыта О.~Виро и представлена в его докладе на конгрессе в
Барселоне \cite{Vir2000}. Исходя из идемпотентного принципа
соответствия, О.~Виро построил кусочно-линейную геометрию
многогранников специального вида в конечномерных эвклидовых
пространствах как результат деквантования Маслова обычной
вещественной алгебраической геометрии. Он указал на важные
приложения к вещественной алгебраической геометрии (например, в
рамках 16-ой проблемы Гильберта о построении вещественного
алгебраического многообразия с предписанными свойствами и
параметрами) и на связь с комплексной алгебраической геометрией и
амебами в смысле И.~М.~Гельфанда, M.~M.~Капранова и А.~В.~
Зелевинского (см. их книгу \cite{GeKaZe94} и статью О.~Виро
\cite{Vir2002}). Затем комплексная алгебраическая геометрия была
деквантована Г.~Михалкиным с тем же результатом; эта новая
``идемпотентная'' (или асимптотическая) геометрия теперь часто
называется {\it тропической алгебраической геометрией},
см., например, \cite{EiKaLi2004,ItKhSh2003,Mi2001,Mi2003,Mi2004,%
Mi2005,RiStTh2005,Shu2002,SpSt2004,St2002}.

Существует естественная связь между деквантованием Маслова и
амебами. Пусть $({\cset}^*)^n$~--- комплексный тор , где ${\cset}^*
= {\cset}\backslash \{0\}$~--- группа ненулевых комплексных чисел по
умножению. Для
 $z = (z_1, \dots, z_n)\in
(\cset^*)^n$ и положительного вещественного числа $h$ обозначим
через $\Log_h(z) = h\log(|z|)$ элемент
\[(h\log |z_1|, h\log |z_2|, \dots,
h\log|z_n|) \in \rset^n.\] Рассмотрим комплексное алгебраическое
многообразие $V\subset (\cset^*)^n$; обозначим через $\maA_h(V)$
множество $\Log_h(V)$. Если $h=1$, то множество $\maA(V) =
\maA_1(V)$ называется {\it амебой} в $V$ в смысле книги
\cite{GeKaZe94},
см. также \cite{Ale2005,Mi2001,Mi2004,Mi2005,PaTs2005,%
Shu2002,So2005,Vir2002}. Амеба $\maA(V)$ является замкнутым
подмножеством в $\rset^n$ с непустым дополнением. Отметим, что
приведенная конструкция зависит от выбора системы координат.

Для простоты выберем в качестве $V$ гиперповерхность
в~$(\cset^*)^n$, определяемую полиномом~$f$; тогда существует
деформация этого полинома $h\mapsto f_h$, порожденная
деквантованием Маслова, и $f_h = f$ для $h = 1$. Пусть $V_h\subset
({\cset}^*)^n$ является множеством нулей полинома $f_h$, и $\maA_h
(V_h) = {\Log}_h (V_h)$. Тогда существует тропическое многообразие
$\mathit{Tro}(V)$, такое что подмножества $\maA_h(V_h)\subset
\rset^n$ сходятся к $\mathit{Tro}(V)$ в метрике Хаусдорфа при
$h\to 0$, см. \cite{Mi2001,Ru2001}. Тропическое многообразие
$\mathit{Tro}(V)$ является результатом деформации амебы $\maA(V)$
и деквантования Маслова многообразия $V$. Множество
$\mathit{Tro}(V)$ получило название {\it скелета} амебы $\maA(V)$.

{\bf Пример \cite{Mi2001}.}  Для прямой  $V = \{\, (x, y)\in
({\cset}^*)^2 \mid x + y + 1 = 0\,\}$ кусочно линейный граф
$\mathit{Tro}(V)$ представляет собой тропическую прямую, см.
рис.~3(a). Амеба $\maA(V)$ показана на рис.~3(b), тогда как
рис.~3(c) демонстрирует соответствующую деформацию амебы.
\begin{figure}[t]
  \label{fig:amoeba}
  \centering
  \includegraphics[width=12cm]{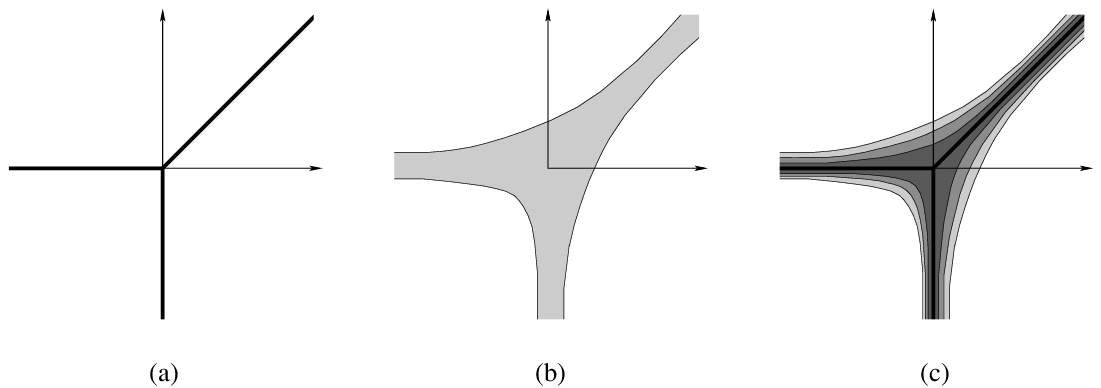}
  \caption{Тропическая прямая и амебы.}
\end{figure}

В важной работе \cite{Kap2000} (см. также
\cite{EiKaLi2004,Mi2001,Mi2004,RiStTh2005}) тропические
многообразия возникают как амебы над неархимедовыми полями.

В 2000 году M.~Концевич заметил, что возможно использование
неархимедовых амеб в перечислительной геометрии, см. [126, раздел
2.4, замечание 4]. Методы тропической геометрии имеют важные
приложения в алгебраической перечислительной геометрии,
инвариантах Громова-Виттена и Вельшенже, см.
\cite{GaMa2005a,ItKhSh2003,%
ItKhSh2003a,ItKhSh2004,ItKhSh2005,Mi2001,Mi2003,Mi2004,Mi2005,Shu2002,Shu2004}.
В частности, Г.~Михалкин представил и доказал в
\cite{Mi2003,Mi2005} формулу, перечисляющую кривые произвольного
рода на торической поверхности (см. также работы
\cite{GaMa2005b,ItKhSh2003,ItKhSh2003a,NiSie2004,Shu2002}).

В последнее время появилось много других работ по тропической
алгебраической геометрии и ее приложениям к обычной (например,
комплексной) алгебраической геометрии, а также к другим областям,
см.,
например \cite{Ale2005,FoGo2003a,FoGo2003b,Izh2005a,NiSie2004,PaSt2003,PaSt2004,%
SzVe2004,Te2004,Vi2004}. Дело в том, что некоторые трудные
традиционные  задачи могут быть сведены к их тропическим версиям,
которые, как можно надеяться, не так трудны.

Заметим, что тропическая геометрия тесно связана с известной
программой M.~Концевича и Я.~Сойбельмана, см., например,
\cite{KonSoi2001,KonSoi2004}.

Введение в тропическую алгебраическую геометрию изложено в работе
\cite{RiStTh2005} (см. также \cite{St2002}). Тем не менее, в целом
сделаны лишь первые шаги в тропической/идемпотентной геометрии, и
задача  систематического построения идемпотентных версий
алгебраической и аналитической геометрии остается открытой.

%% {\bf 12. Correspondence principle for algorithms and their computer
%% implementations.}
\section{Принцип соответствия для алгоритмов и их компьютерная реализация}

В идемпотентной математике имеется большое количество важных прикладных алгоритмов,
см., например, \cite{BaCoOlQu92,Bu2005,Ca79,Cu95,CuMe80,FiRo93,%
FlMc2000,FlMc2005,GoMi79,GoMi2001,Gun98a,ItKhSh2003,KiRo2004,Ki2001,%
Ko2001,KoMa97,LiMa95,LiMa98,LiMaE2000,LiMaRo2000,LiSo2000,LiSo2001,%
LoPe2005,Mi2003,Mi2005,RiStTh2005,Roy92,St2002,Vor63,Vor67,Vor70,%
WoOl2005,ZiK2005,ZiU81}. Идемпотентный принцип соответствия
справедлив как для алгоритмов, так и для их программной и
аппаратной реализации
\cite{LiMa95,LiMa96,LiMa98,LiMaE2000,LiMaRo2000}. В частности,
благодаря принципу суперпозиции, особую важность приобретают
алгоритмы линейной алгебры. Хорошо известно, что алгоритмы
линейной алгебры удобны для параллельных вычислений; их
идемпотентные аналоги также допускают распараллеливание. Поэтому
имеется регулярный способ использовать параллельные вычисления для
решения многих задач, включая основные задачи оптимизации. При
этом удобно использовать универсальные алгоритмы, не зависящие от
конкретного полукольца и от его конкретной компьютерной модели.
Программная реализация полукольцевых универсальных алгоритмов
строится на  объектно-ориентированном и обобщенном (generic)
программировании; программные модули имеют дело с абстрактными
(переменными) операциями и типами данных, см.
\cite{LiMa95,LiMa98,LiMaE2000,LiMaRo2000,LoPe2005}.

Наиболее важные и стандартные алгоритмы имеют много аппаратных
реализаций в виде технических устройств или специальных
процессоров. Такие устройства могут стать прототипами новых
аппаратных устройств, порожденных заменой традиционных
арифметических операций на их полукольцевые аналоги, см.
\cite{LiMa95,LiMa98,LiMaRo2000}. Удачные и эффективные технические
идеи и решения могут быть перенесены из прототипов в новые
аппаратные устройства. Таким образом принцип соответствия
порождает регулярный  эвристический метод для конструирования
многочисленных устройств.

%%{\bf 13. Idempotent interval analysis.}
\section{Идемпотентный интервальный анализ}

Идемпотентная версия традиционного интервального анализа изложена
в  работах \cite{LiSo2000,LiSo2001}. Пусть  $S$~--- идемпотентное
полукольцо, наделенное стандартным порядком. {\it Замкнутым
интервалом} в $S$ назовем множество вида ${\bf x} = [{\bf
\underline x}, {\bf \bar x}] = \{ x\in S\mid {\bf \underline
x}\cle x \cle{\bf \bar x}\}$, где элементы $\bf \underline x\cle
\bf \bar x$ называются {\it нижней} и  {\it верхней границами}
интервала $\bf x$. {\it Слабым интервальным расширением } $I(S)$
полукольца $S$ назовем множество всех замкнутых интервалов $S$ с
операциями $\oplus$ и $\odot$, определенными следующим образом:
$\bf x\oplus\bf y = [\bf \underline x \oplus\bf \underline y, \bf
\bar x\oplus\bf\bar y]$, $\bf x\odot\bf y = [\bf \underline
x\odot\bf \underline y, \bf \bar x\odot\bf \bar y]$. Множество
$I(S)$ является новым идемпотентным полукольцом относительно этих
операций. Доказано, что основные интервальные задачи идемпотентной
линейной алгебры имеют полиномиальную сложность, тогда как в
обычном интервальном анализе задачи такого типа в общем случае
имеют NP-сложность. Точные интервальные решения для дискретного
стационарного уравнения Беллмана (см. матричные уравнения,
рассмотренные выше, в разделе 6) и для соответствующих
оптимизационных задач построены Г.~Л.~Литвиновым и
А.~Н.~Соболевским  \cite{LiSo2000,LiSo2001}. Близкие результаты
были получены в работе \cite{CeCu2002}.

%%{\bf 14. Relations to the KAM theory and optimal transport.}

\section{Связь с КАМ теорией и оптимальным транспортом}

Предмет теории Колмогорова-Арнольда-Мозера (КАМ)  можно
формулировать как изучение инвариантных подмножеств в фазовых
пространствах неинтегрируемых гамильтоновых динамических систем,
имеющих ту же динамику и ту же степень регулярности, что и
интегрируемые системы (квазипериодическое поведение). В последнее
время заметный прогресс был достигнут благодаря вариационному
подходу, в котором динамика задается не гамильтонианом, а
лагранжианом. Соответствующая теория была инициирована работами
С.~Обри (S.~Aubry) и Дж.~Н.~Мезера (J.~N.~Mather) и недавно
названа {\it слабой КАМ теорией} в работах А.~Фати (A.~Fathi)
(готовится к выходу его монография ``Weak KAM Theorems in
Lagrangian Dynamics'' , см., также работы
\cite{KKS2005a,KKS2005b,So99a,So99b}). Минимизация некоторого
функционала вдоль траекторий движущихся частиц является
центральным пунктом другой теории, а именно, теории оптимального
транспорта, которая интенсивно развивается. Эта теория восходит к
работе Г.~Монжа ``О выемках и насыпях'' (1781). Современная версия
теории носит название {\it теория оптимального транспорта
Монжа--Ампера--Канторовича
  (MAK) }
  (название возникло после работы Л.~В.~Кан\-торовича
\cite{K42}).  Имеется сходство между этими двумя теориями и
имеются взаимосвязи с задачами идемпотентного функционального
анализа (например, с задачей о собственных функциях
``идемпотентных''  интегральных операторов, см. \cite{So99a}).
Приложения теории оптимального транспорта к теории обработки
данных в космологии рассматриваются, например, в работах
\cite{BFHLMMS2003,FMMS2002}.

%%{\bf 15. Relations to logic, fuzzy sets, and possibility theory.}

\section{Связь с логикой, нечеткими множествами и~теорией возможностей}

Пусть $S$~--- идемпотентное полукольцо с нейтральными элементами
$\0$ и $\1$ (напомним, что $\0\neq \1$, см.\ разд.~2). Тогда
булева алгебра $\bset = \{\0, \1\}$ является естественным
подполукольцом в $S$. Поэтому $S$ может рассматриваться как
обобщенная (расширенная) логика с логическими операциями $\oplus$
(дизъюнкция) и $\odot$ (конъюнкция). Идеи такого рода обсуждаются
во многих статьях и монографиях, посвященных обобщенным версиям
классической и, особенно, квантовой логики, см., например,
\cite{DiGe2005,Gol99,KlPa2004,Ro90,Ro96}.

Пусть $\Omega$ является так называемым универсумом, состоящим из
``элементарных событий''. Обозначим через ${\CF}(S)$ множество
функций, заданных на $\Omega$ и принимающих значения в $S$; тогда
${\CF}(S)$ представляет собой идемпотентное полукольцо
относительно поточечного сложения и умножения функций. Назовем
элементы из ${\CF}(S)$  {\it обобщенными нечеткими множествами},
см. работы \cite{Gol99,Li2004}. Если $S = {\pset}$, где $\pset$ -
отрезок $[0,1]$ с операциями $\oplus = \max$ и $\odot = \min$, то
мы получим классическое определение нечеткого множества
(Л.~A.~Заде \cite{Za65}). Разумеется, функции из ${\CF}(\pset)$,
принимающие значения в булевой алгебре $\bset = \{0, 1\}\subset
{\pset}$ соответствуют обычным множествам из $\Omega$, а
полукольцевые операции отвечают стандартным операциям над
множествами. В общем случае функции из ${\CF}(S)$, принимающие
значения в $\bset = \{{\0}, {\1}\}\subset S$, могут
рассматриваться как подмножества в $\Omega$. Если $S$~--- решетка
(т.e. $x\odot y = \inf \{x, y\}$ и $x\oplus y = \sup \{x, y\}$),
то обобщенные нечеткие множества совпадают с $L$-нечеткими
множествами в смысле Дж.~А.~Гогена \cite{Go67}. Множество
интервалов $I(S)$ является идемпотентным полукольцом, (см. раздел.
11), поэтому элементы из ${\CF}(I(S))$ могут рассматриваться как
интервальные (обобщенные) нечеткие множества.

Как хорошо известно,  классическая теория нечетких множеств лежит
в основе теории возможностей \cite{Za78,DuPrSa2001}. Разумеется,
исходя из обобщенной теории нечетких множеств, можно аналогичным
образом развить обобщенную теорию возможностей (см., например,
\cite{DuPrSa2001,KlPa2004,Li2004}. Обобщенные теории могут быть
некоммутативными; они выглядят более качественными и менее
количественными по сравнению с обычной теорией, представленной в
работах \cite{Za65,Za78}. Мы видим, что идемпотентный анализ и
теория (обобщенных) нечетких множеств имеют дело с  одними и теми
же объектами, а именно с функциями, принимающими значения в
полукольцах. Тем не менее, основные задачи и методы для обоих
теорий могут отличаться (как это происходит с теорией меры и
теорией вероятностей).

%%{\bf 16. Relations to other areas and miscellaneous applications.}
\section{Приложения к другим областям и к смежным наукам}

Выше обсуждались некоторые приложения идемпотентной математики к
различным теоретическим и прикладным областям математической науки
и взаимосвязи идемпотентной математики с другими областями.
Разумеется, очень естественной областью для приложений идей и
методов идемпотентной математики являются теория оптимизации и
теория оптимального управления. По этой теме имеется  очень
хорошая обзорная статья В.~Н.~Колокольцова \cite{Ko2001};
cм., также работы \cite{BaCoOlQu92,Bu2005,Ca79,CoGaQu99,CoQu94,Cu79,Cu95,%
CuMe80,De97,DeDo98,DeDo2001,Fl2002,Fl2004,FlMc2000,FlMc2005,%
GoMi79,GoMi2001,Gun98a,LiMa95,LiMa98,LiSo2000,LiSo2001,LoQuMa2005,%
Mas86,%
Mas87a,Mas87b,MaKo94,MaSa92,MaVo88,Qu94,Vor63,Vor67,Vor70,%
WoOl2005,ZiK2005,ZiU81}.

Существует много интересных приложений к дифференциальным
уравнениям и к стохастическим дифференциальным уравнениям, см.,
например,
\cite{Fl2002,Fl2004,FlMc2000,FlMc2005,Gun98a,Ko96,Ko2000,Ko2001a,%
KoMa97,Mas86,Mas87a,Mas87b,MaKo94,MaSa92,Pap2005,So99a,So99b}.

Приложения к теории игр обсуждаются, например, в
\cite{KoMal97,KoMa97,MaKo94}. Существуют интересные приложения в
биологии (биоинформатике), см., например,
\cite{FiRo93,PaSt2004,Roy92}. Связь с математической морфологией и
возможные приложения рассматриваются в статье П.~Дел Морала и
М.~Дуази \cite{DeDo2001} и особенно в расширенной препринтной
версии этой работы. Имеется много интересных взаимосвязей с
физикой и соответствующих приложений (квантовая и классическая
физика, статистическая физика, космология и т.д.) см., например, раздел 6 выше и работы \cite{Ko2000,%
KoMa97,LiMaSh2001,LiMaSh2002,LoQuMa2005,Nu91,Qu97,ChDu87}.

Ряд важных взаимосвязей и приложений относится и к чисто
математическим областям. Так называемая тропическая комбинаторика
обсуждается в большой обзорной статье  А.Н.Кириллова
\cite{Ki2001}, см. также \cite{Bu2005,ZiU81}. Интересные
приложения тропических полуколец к традиционной теории
представлений рассматриваются в работах
\cite{BeFoZe96,BeZe2001,Ki2001}. Тропическая математика тесно
связана с весьма популярной теорией кластерных алгебр, основанной
С.~Фоминым и A.~Зелевинским, см. их обзорную работу
\cite{FoZe2004}. В обоих случаях прослеживается связь с
традиционной теорией представлений групп Ли и с близкими задачами.
Имеются важные связи  с выпуклым анализом и с дискретным выпуклым
анализом, см., например, \cite{AkGa2003,CoGaQuSi2005,Cu95,DaKo2004,%
DeSt2003,LiSh2005,MaTi2003,MaSa92,Si97,ZiK77,ZiK79a,ZiK79b}.
Некоторые результаты по сложности идемпотентных и тропических
вычислений могут быть найдены в работах
\cite{KiRo2005a,LiSo2000,LiSo2001,The2004}. Интересные\pagebreak\ приложения
тропической алгебры к теории кос и к преобразованиям Янга-Бакстера
(в смысле \cite{Buh98}) можно найти в работах
\cite{DehDy2002,Dy2002,DyWi2004}.

Начиная с Н.~Н.~Воробьева \cite{Vor63,Vor67,Vor70}, многие авторы
используют, явно или неявно, конструкции и результаты
идемпотентной математики в математической экономике, см.,
например, \cite{DaKoMu2001,KoMal97,MaSa92,ZiK76,ZiU81}.